\documentclass[twoside]{article}

\usepackage[preprint]{aistats2027}

\usepackage{booktabs}
\usepackage{multirow}
\usepackage{adjustbox}
\usepackage{graphicx}
\usepackage{xcolor}
\usepackage{colortbl}
\usepackage{amssymb}
\usepackage{amsmath}
\usepackage{mathrsfs}
\usepackage{amsthm}
\usepackage{hyperref}
\usepackage[normalem]{ulem}
\usepackage{enumitem} 
\usepackage{float}

\usepackage{amsmath}
\usepackage{algorithm}
\usepackage{algorithmic}

\newtheorem{definition}{Definition}
\newtheorem{proposition}{Proposition}
\newtheorem{theorem}{Theorem}

\newtheorem{example}{Example}

\newtheorem*{theorem*}{Theorem}
\newtheorem*{proposition*}{Proposition}
\newtheorem*{definition*}{Definition}
\newtheorem*{example*}{Example}

\usepackage[round]{natbib}

\newcommand{\R}{\mathbb{R}}
\newcommand{\Info}{\mathcal{F}}
\newcommand{\Sig}{S}   
\newcommand{\SH}[1]{S_{\mathrm{H},#1}}    
\newcommand{\SO}[1]{S_{\mathrm{O},#1}} 
\newcommand{\dyn}{\mathrm{dyn}}   
\newcommand{\point}{\mathrm{pt}}  

\newcommand{\hbetaH}[1]{\hat\beta^{(#1)}_{\mathrm{H}}}
\newcommand{\hbetaO}[1]{\hat\beta^{(#1)}_{\mathrm{O}}}
\newcommand{\fH}[2]{f^{(#1)}_{\mathrm{H},#2}}
\newcommand{\fO}[2]{f^{(#1)}_{\mathrm{O},#2}}

\usepackage{amsmath}

\DeclareMathOperator*{\argmin}{arg\,min}

\begin{document}

%

%

\twocolumn[

\aistatstitle{Encoding Time Series Dynamics with Sliding-Window Signatures for Online Forecasting}

\aistatsauthor{ Nina Drobac \And Margaux Br\'eg\`ere \And  Joseph de\ Vilmarest \And Olivier Wintenberger}

\aistatsaddress{Sorbonne Université \And Sorbonne Université \\ EDF R\&D  \And Viking Conseil \And Sorbonne Université} ]

\begin{abstract}
    Accurate forecasting often relies on both past target values and exogenous variables as inputs. In many applications, however, the time series of interest are non-stationary and the relationship between target and inputs evolves over time, requiring forecasters that adapt online. Moreover, inputs affect the target not only through their current values but also through their recent joint dynamics. We address these challenges in three steps. First, to encode the series dynamics, we rely on signatures of the input paths over sliding windows. Motivated by the universality of linear functionals of signatures, we fit a regularized linear regression on these features to capture nonlinear, path-dependent effects. Second, the simplicity of the linear model allows us to update it online as new data arrives, using a Kalman filter to track dynamic effects. Third, to capture dependencies at several time scales, we combine forecasters built on different window lengths through online aggregation. For French national electricity load forecasting, our method outperforms various state-of-the-art approaches. It also remains competitive across several benchmark datasets with exogenous variables.
\end{abstract}

\section{INTRODUCTION}

Many sectors of human activity rely on accurate forecasts of a target of interest
whose future values depend not only on its own past, but also on exogenous
covariates. Recent neural forecasting models report substantial gains in accuracy when incorporating informative exogenous variables \citep{olivares2023nbeatsx,das2024tide,wang2024timexer}. The same holds for
pretrained foundation models: Chronos-2 obtains its largest improvements over
competing models on covariate-informed tasks \citep{ansari2025chronos2, shchur2026fevbench}. Alongside these recent approaches, a well-established way to forecast with exogenous
information is regression, such as generalized additive models (GAMs), which express the target as a sum of smooth functions of the covariates \citep{wood2017generalized}.

\textit{Covariate dynamics.} Beyond their current values, the dynamics of the covariates also matter: a slow decrease of a covariate can affect the target differently from a sudden drop ending at the same value. In electricity load forecasting, for instance, the thermal inertia of buildings makes demand respond to the recent temperature trajectory rather than to the current temperature alone \citep{pierrot2011short}. In all the approaches above, however, covariates mostly enter pointwise: their dynamics must either be engineered by hand, e.g.\ through exponentially smoothed temperatures \citep{goude2013local}, or be learned implicitly from large amounts of data.

\textit{Dynamics evolving over time.} Real-world time series are rarely stationary, and the relationship between the target and the covariates 
can itself evolve over time. Such drifts call for models that are updated as new data arrive. Neural networks find this online setting challenging, as they are costly to retrain, and incremental updates expose them to catastrophic forgetting \citep{pham2023learning,vandeven2025continual}. Foundation models, in turn, adapt only through their context window or fine-tuning. Simpler models, by contrast, can be updated sequentially at low cost: the coefficients of linear models and GAMs can be tracked by Kalman filtering \citep{obst2021adaptive,devilmarest2022state}, and several models can be combined by online aggregation of experts \citep{cesabianchi2006prediction,gaillard2015forecasting}.

\textit{Path signatures.} 
Simpler models, such as linear ones, are therefore better suited to the online
setting, provided that they are given rich enough features that encode the dynamics of the
covariates series. Path signatures provide such a feature set.  Introduced by \citet{Chen_1958} and developed by \citet{Lyons_2007} as a cornerstone of rough path theory, the signature is a representation of a path, defined as the collection of its iterated integrals. It enjoys a universal approximation property: any continuous function of a path can be approximated by a linear function of its signature. Interpreting observations as paths and computing their signatures has therefore become an effective way to obtain non-parametric features for sequential data \citep{chevyrev2026primer,fermanian2021thesis}. 
In forecasting, they have been combined with transformers \citep{morenopino2025roughtransformerslightweightcontinuous} and used as inputs to linear models, for instance to nowcast macroeconomic series \citep{cohen2023nowcasting} and to forecast freight transportation rates at Amazon \citep{gu2025transportation}.
For a time series observed
sequentially, one must choose the interval over which the signature is computed. Computing it over the whole past weighs all observations equally and turns a stationary time series into non-stationary features. Two solutions have been proposed: exponentially fading memory signatures discount the past of the path \citep{abijaber2025efm}, while sliding-window signatures discard observations outside a window of fixed length \citep{drobac2025sliding}. We build on the latter, as they are computed as standard signatures and their windows have a direct interpretation as a time scale; extending the approach to the non-stationary setting.


\paragraph{Our contributions.} 

\textit{A novel forecasting framework.} We propose SigAdapt, a lightweight and theoretically grounded forecasting model that addresses three difficulties of forecasting with covariates.
(i)~Designing time series features is challenging when the series dynamics change over time. SigAdapt relies on sliding-window signatures as an expressive feature set that naturally encodes these dynamics. Furthermore, since linear functionals of the signature are universal approximators, a linear
model is a natural choice on these features.
(ii)~The relationship between the covariates and the target drifts over time. The simplicity of the linear model allows a Kalman filter to track the dynamic effects of the signature features online.
(iii)~The relevant memory length is unknown and may change over time. SigAdapt
combines models built on different window lengths by online aggregation of experts.

\textit{Cheap online updates.} Each step of SigAdapt is updated sequentially: the signature features are computed efficiently as the window slides thanks to algebraic properties, the Kalman filter operates on a low-dimensional summary of the signature effects, and the online aggregation methods are inherently fast.

\textit{Empirical evaluation.} On French national electricity load and on benchmark electricity datasets, SigAdapt outperforms neural, foundation-model and additive model competitors when the covariates are informative for the target. Furthermore, an ablation study illustrates that each SigAdapt component contributes to its performance.

\section{Signatures}

In this section, we introduce signatures and frame them as features for time series. For a more formal treatment of signatures, see Appendix \ref{app:sig_theory}.

\subsection{Path signatures}
Signatures are defined for continuous paths, i.e. continuous mappings from $[s,t]$ to $\mathbb{R}^d$. Let $d \geq 1$ and let $\mathrm{BV}([s,t],\mathbb{R}^d)$ denote the set of continuous paths $x : [s,t] \to \mathbb{R}^d$ with finite $1$-variation
\begin{equation*}
    \|x\|_{1\text{-}var,[s,t]} = \sup_{\mathcal{D}} \sum_{j=0}^{m-1} \big\| x_{u_{j+1}} - x_{u_j} \big\|,
\end{equation*}
where the supremum is taken over all time partitions $\mathcal{D} = \{s = u_0 < \cdots < u_m = t\}$ of $[s,t]$ and $\|x\|^2=\sum_{i=1}^d(x^{(i)})^2$, $x\in \R^d$. For such paths, the iterated integrals below are well-defined.

\begin{definition}[Signature]
\label{def:signature}
    Let $x \in \mathrm{BV}([s,t],\mathbb{R}^d)$. For a multi-index $I = (i_1, \ldots, i_k) \in \{1, \ldots, d\}^k$ of length $|I| = k$, the signature coefficient of $x$ on $[s,t]$ is 
    \begin{equation*}
        S^{I}\left( x_{[s,t]} \right) = \idotsint \displaylimits_{s < u_1 < \cdots < u_k < t} \mathrm{d}x^{(i_1)}_{u_1} \cdots \mathrm{d}x^{(i_k)}_{u_k},
    \end{equation*}
    with $S^{\emptyset} = 1$. The signature of level $k$ is   $S^{(k)}(x_{[s,t]}) = \big( S^I(x_{[s,t]}) \big)_{|I| = k} \in \mathbb{R}^{d^k}$. The collection of all levels defines the path signature $S(x_{[s,t]}) = \big( S^{(k)}(x_{[s,t]}) \big)_{k \geq 0}$, and the truncated signature of order $N$ is
        \begin{equation*}
        S^{\leq N}\left( x_{[s,t]} \right) = \Big( 1, S^{(1)}\left( x_{[s,t]} \right), \ldots, S^{(N)}\left( x_{[s,t]} \right) \Big) \in \mathbb{R}^{s_d(N)},
    \end{equation*}
    where $s_d(N) = \sum_{k=0}^{N} d^k = \frac{d^{N+1}-1}{d-1}$.
\end{definition}

Two results make signatures useful as features. First, their levels decay factorially in norm, which justifies truncation in practice.

\begin{proposition}[Factorial decay, \citealp{lyons2014rough}]
\label{prop:decay}
    Let $x \in \mathrm{BV}([s,t],\R^d)$. Then, for every $k \geq 0$,
    $\big\| S^{(k)}\left( x_{[s,t]} \right) \big\| \leq \|x\|_{1\text{-}var,[s,t]}^k / k!$.
\end{proposition}

Second, any continuous function of a path can be uniformly approximated by a linear 
function of its truncated signature, which turns a nonlinear learning problem into a 
linear one.

\begin{theorem}[Universal approximation]
\label{thm:uat}
    Let $K$ be a compact subset of $\mathrm{BV}([s,t],\R^d)$, equipped with the
    $1$-variation topology, such that every $x \in K$ satisfies $x_s = 0$ and
    $x^{(1)}_u = u - s$ for all $u \in [s,t]$. Let $f : K \to \R$ be continuous.
    Then, for every $\varepsilon > 0$, there exist $N \geq 1$ and
    $\theta \in \R^{s_d(N)}$ such that
    \[
        \sup_{x \in K} \big| f(x) - \langle \theta, S^{\leq N}(x_{[s,t]}) \rangle \big|
        \leq \varepsilon .
    \]
\end{theorem}

\subsection{Signatures for time series}
\label{sec:sigwindows}

To use signatures as features for a discrete-time series, we first need to define a path this series traces, and we do so by linearly interpolating the observations. However, the signature is 
invariant under translations $x \mapsto x + c$ and under reparametrisations of $[s,t]$, meaning it does not encode information on the first point (level) of the path nor on the speed at which it is traversed (Proposition \ref{prop:sig_properties} Appendix~\ref{app:sig_theory}). We restore both with two standard augmentations 
\citep{morrill2021generalised}. The time augmentation adds time as the first coordinate, and the basepoint augmentation adds the origin as the first point. Together, they ensure that the signature characterises the embedded observations \citep{hambly2010uniqueness}.

Formally, let $(X_u)_{u \in \mathbb{Z}}$ be a discrete-time series in $\R^{d}$. For integers $s < t$, let $v = t - s$. With a slight abuse of notation, we denote by $X_{[s,t]} : [0,1] \to \R^{d+1}$ the time- and basepoint-augmented path obtained by linearly 
interpolating the points
\begin{equation*}
    x_0 = 0, \quad x_{j/(v+1)} = \Big( \tfrac{j}{v+1},\ X_{s+j-1} \Big), \quad j = 1, \ldots, v+1,
\end{equation*}
and we write $S(X_{[s,t]})$ for its signature on $[0,1]$. By construction, paths defined this way meet the assumptions required in Theorem~\ref{thm:uat}.

At each time $t$, a common choice for the time series embedding is the \emph{expanding-window signature} $S(X_{[t_0,t]})$, computed from a fixed initial time $t_0$, to which new observations are appended as $t$ increases. However, it weighs all past observations equally and does not discount the distant past. Instead, to focus on recent and relevant observations, we use the \emph{sliding-window signature} $S(X_{[t-w,t]})$ of window length $w \in \mathbb{N}$: as $t$ 
increases, the newest observation enters the window and the oldest one leaves it, so that the signature only encodes the recent past.

\section{Methodology}
\label{sec:method}

\begin{figure*}[t]
\centering
\includegraphics[width=0.8\textwidth]{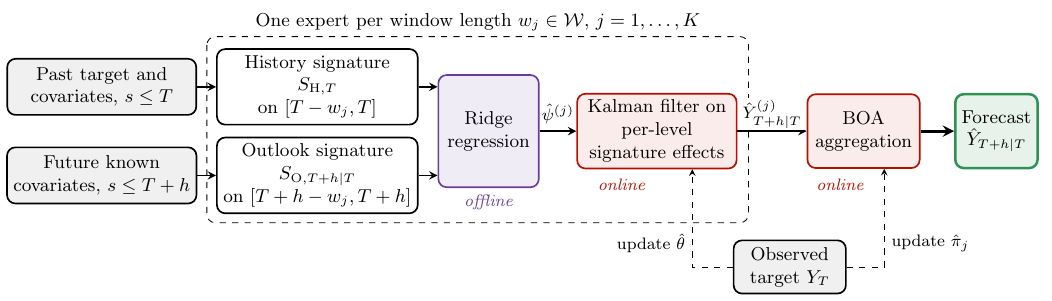}
\caption{Flowchart of SigAdapt at forecast origin $T$, to obtain a prediction $t = T+h$.
For each of the $K$ window lengths $w_j$, the history and outlook signatures serve as series features in a ridge regression fitted once on the training period (purple). Per-level signature effects are adapted online by a Kalman filter and the $K$ resulting forecasts are
combined by Bernstein Online Aggregation (red).}
\label{fig:sigadapt_flowchart}
\end{figure*}

We now present SigAdapt, a forecasting methodology composed of four steps. First, sliding-window
signatures encode the recent trajectories of the target and the covariates into informative features (Section~\ref{sec:encoding_sig}). These features then enter a linear regression model, fitted offline on the training
period to obtain coefficient estimates (Section~\ref{sec:lin_model}). A Kalman filter on signature-derived effects lets the linear model adapt to changes in the relationship between target and features (Section~\ref{sec:kalman}). Finally, online
aggregation combines models built on different window lengths (Section~\ref{sec:aggregation}). As this section is notation-heavy, we refer to Appendix~\ref{app:glossary} for the glossary. The whole methodology is summarized in Figure~\ref{fig:sigadapt_flowchart}.

\subsection{Forecasting setting}

We place ourselves in the setting of univariate forecasting with exogenous covariates for discrete-time time series. Let $(Y_s)_{s \in \mathbb{Z}}$ denote the target time series, with $Y_s \in \R$. The covariates can include both continuous and discrete random variables, and we distinguish two types according to when they are available: \emph{past covariates} $(X_s)_{s \in \mathbb{Z}}$, with $X_s \in \R^{d_X}$, which are observed only up to the forecast origin, and \emph{future known covariates} $(Z_s)_{s \in \mathbb{Z}}$, with $Z_s \in \R^{d_Z}$, whose values are known in advance. Examples of future covariates are calendar variables, such as the day of the week, and forecasts issued at the origin, such as temperature forecasts. At a forecast origin $T$, we forecast the target at time $t = T + h$, for a lead time $h \geq 1$. The information available at the origin is given through the $\sigma$-algebra
\begin{equation}\label{eq:infoset}
  \Info_T = \sigma\big( Y_s, X_s : s \le T;\ Z_s : s \le T + h \big).
\end{equation}
A forecast $\hat Y_{t\mid T}$ is any $\Info_T$-measurable random variable.

\subsection{Encoding the dynamics with signatures}\label{sec:encoding_sig}

Sliding-window signatures are designed to provide an informative summary of a path
inside a window, capturing its evolution rather than its value at a single time. It is therefore the covariates whose recent trajectory influences the target that benefit from the signature transformation. When only the current value of a covariate is informative, it can enter the model directly. With this in mind, we split each covariate vector into two sub-vectors: \emph{dynamic} covariates, $X^{\dyn}_s$ and $Z^{\dyn}_s$, which enter a signature path, and \emph{pointwise} covariates, $X^{\point}_s$ and $Z^{\point}_s$, which enter the model directly through their value at a single time.

Similar to recent covariate-aware neural forecasters, which combine a look-back window (context) with the known future covariates up to the target time \citep{ansari2025chronos2}, we compute two signatures for each forecast, each on its own time window but with the same window length $w$ for simplicity. The \emph{history signature} encodes the recent trajectory of the target and the dynamic covariates on the window $[T - w,\, T]$, while the \emph{outlook signature} describes how the dynamic future known covariates evolve on the window $[t  - w,\, t]$ ending at the target time $t$ :
\begin{align}
  \SH{T} &= \Sig^{\le N_{\mathrm H}}\bigl(H_{[T - w,\, T]}\bigr), \label{eq:hist}\\
  \SO{t\mid T} &= \Sig^{\le N_{\mathrm O}}\bigl(O_{[t - w,\, t]}\bigr), \label{eq:outlook}
\end{align}
where $H_s = \bigl(Y_s,\, X^{\dyn}_s,\, Z^{\dyn}_s\bigr)$ and $O_s = Z^{\dyn}_s$, $s \in \mathbb{Z}$, and $H_{[T-w,\,T]}$ and $O_{[t-w,\,t]}$ are the corresponding augmented paths of Section~\ref{sec:sigwindows}. A natural extension is to allow distinct lengths $w_{\mathrm{H}}$ and $w_{\mathrm{O}}$.

\begin{figure}[H]
\centering
\includegraphics[width=\linewidth]{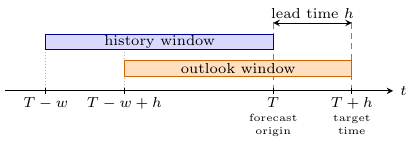}
\caption{Windows used at origin $T$ for lead time $h$.} 
\label{fig:windows}
\end{figure}

We illustrate the two signature windows in Figure \ref{fig:windows}. The history signature captures not only the autoregressive dynamics of the target, but also its relationship with the covariates: its higher-level terms, such as the level-$2$ signed areas, describe the joint evolution of the target and each covariate. The outlook signature complements it with the evolution of the future known covariates up to the target time. Both signatures are therefore $\Info_T$-measurable and, together with the pointwise covariates, form the features of the linear model introduced in the next section.

\subsection{Linear model}\label{sec:lin_model}

The use of signature features is motivated by the universal approximation (Theorem \ref{thm:uat}): instead of learning a complex nonlinear function of the path, we can approximate it by a linear functional of the truncated signature. 

We therefore consider static forecasts of the form
\begin{equation}\label{eq:linmodel}
  \hat Y_{t \mid T}
   = \alpha + \big\langle\gamma, U_{t\mid T}\big\rangle + \big\langle \beta^{\mathrm H}, \SH{T} \big\rangle + \big\langle \beta^{\mathrm O}, \SO{t\mid T} \big\rangle,
\end{equation}
with $U_{t\mid T} = \bigl(X^{\point}_T,\, Z^{\point}_t\bigr)$ denoting the pointwise covariates.

The parameters $\psi = (\alpha, \gamma, \beta^{\mathrm H}, \beta^{\mathrm O})$ are estimated once, on a training period. 
Signature features are strongly collinear, which makes the design matrix ill-conditioned. We therefore use ridge regression. Since signature coefficients at different levels live on different scales (Proposition~\ref{prop:decay}), we first standardize each feature by its mean and standard deviation over the training period, so that the penalty treats all features alike. With $\mathcal{T}_{\mathrm{tr}}$ the set of training forecast origins, the estimator is given as
\begin{equation}\label{eq:ridge}
  \hat\psi = \argmin_{\psi}\
  \sum_{T\in\mathcal{T}_{\mathrm{tr}}} \bigl(Y_{T+h} - \big\langle\psi, \varphi_{T+h\mid T}\big\rangle\bigr)^2 + \lambda \|\psi\|_2^2,
\end{equation}
with $\varphi_{t\mid T} = \bigl(1,\, U_{t\mid T},\, \SH{T},\, \SO{t\mid T}\bigr)$.  The regularization parameter $\lambda$ is selected on a validation period. 

Beyond being compact and interpretable, the linear model extends naturally to the time-varying setting of Section~\ref{sec:kalman}, since it can be adapted as new observations arrive, at a low computational cost. The ridge estimate~\eqref{eq:ridge} initializes the adaptive model presented in the following section.

\subsection{Kalman filter}\label{sec:kalman}

In many real-world settings, the relationship between the target and the inputs changes over time, so the coefficients estimated on the training set can become outdated. To update our model coefficients as new data arrive, we let the coefficients evolve as a latent state vector, and we estimate its {\it a posteriori} distribution. We rely on a linear Gaussian state-space model where both the state innovation and the conditional distribution of the observation given the state vector are Gaussian. The state posterior distribution is then given by the Kalman filter \citep{kalman_new_1961}.


Following \citet{obst2021adaptive}, who apply the Kalman filter to aggregated effects of a generalized additive model rather than to each of its parameters, we filter one effect per signature level, which respects the natural ordering of the signature. Let $\SH{T}^{(k)}$ and $\SO{t\mid T}^{(\ell)}$ denote the standardized level-$k$ and level-$\ell$ coordinates of the history and outlook signatures in~\eqref{eq:ridge}, and $\hat\beta^{\mathrm H}_k$ and $\hat\beta^{\mathrm O}_\ell$ the corresponding ridge coefficients. The per-level \emph{signature dynamic effects} are
\begin{equation*}
\begin{aligned}
\fH{k}{T} &= \big\langle \hbetaH{k},\, \SH{T}^{(k)} \big\rangle, && 1 \leq k \leq N_\mathrm{H},\\
\fO{\ell}{t\mid T} &= \big\langle \hbetaO{\ell},\, \SO{t\mid T}^{(\ell)} \big\rangle, && 1 \leq \ell \leq N_\mathrm{O},
\end{aligned}
\end{equation*}
so that $\sum_k \fH{k}{T} = \langle \hat\beta_{\mathrm{H}}, \SH{T} \rangle$ and 
$\sum_\ell \fO{\ell}{t\mid T} = \langle \hat\beta_{\mathrm{O}}, \SO{t\mid T} \rangle$.

We split the pointwise covariates into $U_{t\mid T} = \big(U^{\mathrm{fix}}_{t\mid T}, U^{\mathrm{ad}}_{t\mid T}\big)$, with ridge coefficients $\hat\gamma = (\hat\gamma_{\mathrm{fix}}, \hat\gamma_{\mathrm{ad}})$. The coefficients $\hat\gamma_{\mathrm{fix}}$ are kept at their batch value, typically for covariates whose effect is not expected to drift while those of $U^{\mathrm{ad}}_{t\mid T} \in \R^{p}$ are adapted, together with one weight per signature level. With the regressor
\begin{equation*}
\begin{aligned}
\xi_{t\mid T} = \big(&1,\ U^{\mathrm{ad}\top}_{t\mid T},\ \fH{1}{T}, \dots, \fH{N_\mathrm{H}}{T},\\
&\fO{1}{t\mid T}, \dots, \fO{N_{\mathrm{O}}}{t\mid T}\big)^\top \in \R^{1+p+N_\mathrm{H}+N_\mathrm{O}},
\end{aligned}
\end{equation*}
we consider the state-space model
\begin{equation}\label{eq:ssm}
\begin{aligned}
\theta_t &= \theta_{t-1} + \eta_t, && \eta_t \sim \mathcal{N}(0, Q),\\
Y_t &= \big\langle\hat\gamma_{\mathrm{fix}} , U^{\mathrm{fix}}_{t\mid T}\big\rangle + \big\langle\theta_t , \xi_{t\mid T}\big\rangle + \varepsilon_t, && \varepsilon_t \sim \mathcal{N}(0, \sigma^2),
\end{aligned}
\end{equation}
where the state $\theta_t$ is indexed by the target time, and $Q$ is diagonal. We fix $\theta_1 \sim \mathcal{N}(\hat\theta_1, P_1)$, with $\hat\theta_1 = (\hat\alpha, \hat\gamma_{\mathrm{ad}}^\top, 1, \dots, 1)^\top$, so that a constant state equal to $\hat\theta_1$ reproduces the batch forecast~\eqref{eq:linmodel}. The resulting online forecast is
\begin{equation}\label{eq:kf-forecast}
\hat Y_{t\mid T} = \big\langle\hat\gamma_{\mathrm{fix}},  U^{\mathrm{fix}}_{t\mid T}\big\rangle + \big\langle\hat\theta_{T\mid T}, \xi_{t\mid T}\big\rangle,
\end{equation}
where $\hat\theta_{T\mid T}$ is the Kalman filter estimate of the state given information up to time $T$. The hyperparameters $(P_1, Q, \sigma^2)$ are selected on the training set by an iterative grid search \citep[Chapter 5]{these_joseph}.

\subsection{Aggregation over window experts}
\label{sec:aggregation}

Choosing the window length $w$ of the history and outlook signatures imposes a
trade-off between reactivity and stability, given that short windows adapt quickly
to changes in the dynamics of the series, whereas long windows yield a more stable
summary of slowly varying effects. Furthermore, as the dynamics of the target and input series change, the best window may vary over time. Rather
than selecting a single $w$ on a validation period, we run the procedure of
Sections~\ref{sec:encoding_sig}--\ref{sec:kalman} in parallel for each window length on a finite
grid $\mathcal{W} = \{w_1, \dots, w_K\}$, and combine the resulting forecasts,
viewed as experts, by expert aggregation.

Let $\hat{Y}^{(j)}_{t\mid T}$ denote the forecast~\eqref{eq:kf-forecast} of $Y_t$
issued at origin $T = t - h$ by the model with window length $w_j$. The aggregated
forecast is the convex combination
\begin{equation}\label{eq:aggregation}
  \hat{Y}_{t\mid T} = \sum_{j=1}^{K} \hat\pi_{j,T}\, \hat{Y}^{(j)}_{t\mid T},
\end{equation}
where the weights $\hat\pi_T = (\hat\pi_{1,T}, \dots, \hat\pi_{K,T}) \in \Delta_K
= \{q \in [0,1]^K : \sum_{j} q_j = 1\}$ are computed from the past forecasts
$\hat{Y}^{(j)}_{s\mid s-h}$ and observations $Y_s$, $s \leq T$, that is, from the
information available at the forecast origin. The weights are computed with
Bernstein Online Aggregation (BOA) \citep{wintenberger2017optimal}, as implemented
in the package \texttt{opera} \citep{gaillard2016opera}. BOA guarantees that, under the
squared loss, the aggregated forecast performs asymptotically as well as the best
fixed convex combination of the experts chosen in hindsight, and in particular as
well as the best single window length.


\subsection{Efficient online updates}\label{sec:online_updates}

\textbf{Signature updates.} The online updates of the Kalman filter and BOA are only fast if the signature features are computed at the same pace. Since consecutive sliding windows have a large overlap, we leverage the algebraic properties of signatures to update them sequentially rather than recomputing them, while still respecting the basepoint augmentation (Algorithm~\ref{algo:update_sig}, Appendix~\ref{app:sig_algo}). For a $d$-dimensional path
truncated at order $N$, an update costs $O(Nd^{N})$ operations, against $O(wNd^{N})$ from scratch, which makes computation on long windows feasible.

\textbf{Complexity.} A Kalman step costs $O(D^{2})$, where
$D = 1 + p + N_{\mathrm H} + N_{\mathrm O}$ is the state dimension, independently of
$d$, and a BOA step costs $O(K)$. One online step thus
costs $O\big(K(Nd^{N} + D^{2})\big)$, independently of $w$. In practice, $d^{N}$ remains moderate: the paths are low-dimensional and, owing to the factorial decay of the signature (Proposition~\ref{prop:decay}), the truncation orders are often kept small (between $4$ and $6$ in our experiments), as is common in literature \citep{fermanian2022functional}.

\section{APPLICATION}
\label{sec:experiments}

We evaluate SigAdapt in two complementary settings. Section~\ref{sec:french_load} is an in-depth case study of short-term forecasting of the French national load, in which all models receive the same hand-crafted calendar covariates, including public and school holidays. This setting allows a detailed analysis of the results, with significance tests and an ablation study of the components of SigAdapt. Section~\ref{sec:bench_data} then assesses the performance of the models when applied directly to public benchmark data with minimal dataset-specific engineering.

The data and code used for this section are publicly available in the GitHub repository: \url{https://github.com/ninadrobac/sigadapt}.

\subsection{Case study: French national load}
\label{sec:french_load}

\paragraph{Data.}
We use the half-hourly French national electricity load (MW) from Eco2mix\footnote{https://www.rte-france.com/eco2mix}, the open data platform of the French transmission system operator RTE. We pair it with temperature (\textdegree C) recorded by Météo-France\footnote{\url{https://donneespubliques.meteofrance.fr}}: we average the 3-hourly observations across French weather stations and interpolate them linearly to match the load grid. The dataset further contains calendar covariates, namely those relating to public holidays. The data cover the period of 1 January 2018 to 30 June 2026.

\paragraph{Forecasting task and setup.}
We consider two-day-ahead forecasts. The forecast of the load at time $T+2\text{ days}$ uses load observations up to the forecast origin $T$, and calendar and temperature values up to the time of the forecast $T+2\text{ days}$. As temperature at the forecasted time is an important predictor of electricity load, we treat temperature as known in the future and use its observed value in place of the meteorological forecasts. In terms of calendar variables, we consider the position in the year, the day of the week, working days, bridges (working days between two non-working ones), the Christmas period, and indicators for the Christmas, July and August school holidays. Hyperparameters are selected on 2022 - 2023 with models fitted on 2018 - 2021; all models are then refitted on 2018 - 2023 and tested from 1 January 2024 to 30 June 2026. We refer to Appendix~\ref{app:french_load} for error metrics definitions and other details.

\paragraph{SigAdapt.}
We fit 48 separate models, one per half-hour of the day. This choice is common in load forecasting with additive models \citep{fan2010short, bregere2022online}, and we further motivate it by three reasons. First, the response of load to the inputs depends on the time of day. Second, each Kalman filter then runs on a daily series and tracks day-to-day drift, whereas a single filter over all half-hours would have to handle the daily cycle. Third, the informativeness of short and long windows varies with the time of day, so the expert weights are also learned per half-hour.
In the notation of Section~\ref{sec:method}, the target $Y_s$ is the load, $Z^{\dyn}_s$ is the temperature, and $Z^{\point}_s$ are the calendar variables; there are no past covariates $X_s$. Hence the history signature 
is computed on the path of load and temperature, and the outlook signature on the path of temperature only, both truncated at depth $6$, while the calendar variables enter the ridge regression directly. The Kalman filter tracks the intercept and one dynamic effect per signature level. The filter is not updated on holidays and bridges, whose residuals do not reflect the level of normal days. The experts based on windows of $w \in \{2, 7, 14, 28, 84, 168\}$ days are aggregated with BOA to obtain the final forecasts.


\paragraph{Competitors.}
We compare SigAdapt with three  competitors, covering different approaches to load forecasting. The first is generalized additive models (GAMs), which are widely used for this task \citep{antoniadis2024statistical}, with smoothed temperatures and lagged loads as hand-crafted features. We report both a static GAM and GAMAdapt, in which, similarly to SigAdapt, a Kalman filter adapts the GAM terms online~\citep{obst2021adaptive}. The other two competitors are neural networks. The first is TiDE~\citep{das2024tide}, a MLP-based encoder-decoder forecaster that supports exogenous covariates and performs well on electricity load \citep{drobac2025sliding}. The second is Chronos-2~\citep{ansari2025chronos2}, a large pretrained foundation model that we apply zero-shot. For both, the context length is selected on the validation set, after which TiDE is refitted on the joint training and validation set. TiDE is the only model whose training is stochastic, so we report the mean and standard deviation of its results over five seeds. All models receive the same calendar covariates, except Chronos-2, which receives only the holiday indicators and infers the daily and weekly patterns from its context. Further details on the implementation can be found in Appendix \ref{app:french_load_models}.

\begin{figure}[t]
\centering
\includegraphics[width=\linewidth]{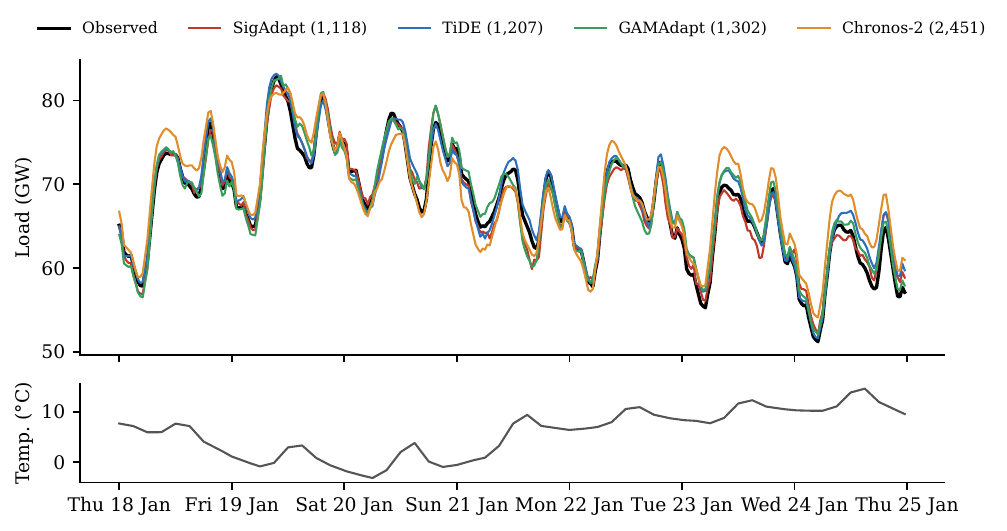}
\caption{Load and forecasts during 18-24 January 2024, the test week with the largest temperature range (lower panel). RMSE in parentheses; TiDE: best seed.}
\label{fig:example_week}
\end{figure}

\paragraph{Results.}
SigAdapt achieves the lowest error on all three metrics (Table~\ref{tab:rte_results}), followed by TiDE, GAMAdapt, the static GAM and Chronos-2. We assess the significance of these differences with one-sided Diebold-Mariano tests on squared errors. SigAdapt is significantly more accurate than GAMAdapt ($p < 0.001$), and so is TiDE ($p = 0.035$). SigAdapt has lower errors than TiDE for all five seeds, but this difference is not significant ($p = 0.21$). Therefore, a lightweight linear model on signature features, updated online, performs on 
par with a deep neural forecaster trained on the same data. Chronos-2, used zero-shot, is significantly less accurate than all other models, including the static GAM ($p < 0.001$). Its errors are largest on public holidays and during cold spells, suggesting that on this long series with strong calendar and weather effects, models fitted to the series keep a clear advantage over 
a zero-shot foundation model.
Figure~\ref{fig:example_week} shows the forecasts during the test week with the largest temperature range: SigAdapt's forecast is the closest to the load curve, which indicates that it captures well the response of the load to the temperature dynamics. Online adaptation matters for both adaptive models. GAMAdapt improves on the static GAM by 15\% 
in RMSE, and the aggregated SigAdapt forecast improves on its best single window per half-hour chosen in hindsight by 7\% (Table~\ref{tab:rte_ablation}).

\begin{table}[ht]
\centering
\resizebox{0.48\textwidth}{!}{%
    \begin{tabular}{lccc}
    \toprule
    Model & RMSE (MW) & MAE (MW) & MAPE (\%) \\
    \midrule
    GAM & 1\,504 & 1\,111 & 2.16 \\
    GAMAdapt & 1\,280 & 939 & 1.85 \\
    SigAdapt & \cellcolor{lightgray} \textbf{1\,168} & \cellcolor{lightgray} \textbf{845} & \cellcolor{lightgray} \textbf{1.66} \\ 
    TiDE & 1\,205 $\pm$ 8 & 883 $\pm$ 9 & 1.72 $\pm$ 0.02 \\
    Chronos-2 & 1\,774 & 1\,256 & 2.45 \\
    \bottomrule
\end{tabular}%
}
\caption{Test RMSE, MAE and MAPE on French national load, on the test period (2024 - 2026). TiDE: mean $\pm$ standard deviation over five seeds.} 
\label{tab:rte_results} 
\end{table}

\paragraph{Ablation study.}
\paragraph{Ablation study.}
To isolate the contribution of each of the three components of SigAdapt, namely the signature 
features, the Kalman filter and the aggregation, we remove one at a time and keep the rest of the 
pipeline unchanged (Table~\ref{tab:rte_ablation}; details in Appendix~\ref{app:rte_ablation}). 
In the model without signatures, the signature features are replaced by the raw inputs: a linear model, adapted by the Kalman filter, uses the temperature at the target time and the load at the forecast origin, two days earlier. To mirror the windows of SigAdapt, we also build one such linear expert per window length $w$, which additionally contains the temperature $w$ days before the target time, and aggregate them. Without the Kalman filter, the static signature experts are aggregated directly. Without aggregation, a single signature window is used for each half-hour, chosen with hindsight as the best on the test period. Removing the signatures increases the errors by about $40\%$, and aggregating the linear experts brings almost no gain: they make nearly identical errors (mean pairwise correlation $0.99$), whereas the signature experts do not ($0.61$), as each window summarizes the temperature and load paths at its own time scale. Removing the Kalman filter or the aggregation increases the errors by $7$--$10\%$, so both forms of online adaptation are needed and complement each other. Finally, without aggregation, the signature model is on par with GAMAdapt, which relies on hand-crafted features to model the influence of past temperature and load.

\begin{table}[ht]
\centering
\resizebox{0.48\textwidth}{!}{%
\begin{tabular}{lcc}
\toprule
\textbf{Variant} & \textbf{RMSE (MW)} & \textbf{MAE (MW)} \\
\midrule
SigAdapt & \textbf{1\,168} & \textbf{845} \\
No signatures: linear model & 1\,625 \,(+39.1) & 1\,188 \,(+40.5) \\
No signatures: aggr.\ linear experts & 1\,607 \,(+37.6) & 1\,178 \,(+39.3) \\
No Kalman filter & 1\,260 \,(+7.9) & 902 \,(+6.7) \\
No aggregation & 1\,257 \,(+7.6) & 926 \,(+9.5) \\
\midrule
GAMAdapt & 1\,280 \,(+9.6) & 939 \,(+11.0) \\
\bottomrule
\end{tabular}%
}
\caption{Ablation on French load, test 2024--2026. In parentheses: increase in error relative to SigAdapt (\%). MAPE and further variants in Appendix~\ref{app:rte_ablation}.}
\label{tab:rte_ablation}
\end{table}

\subsection{Evaluation on public benchmarks}
\label{sec:bench_data}

\paragraph{Data.} We further evaluate SigAdapt on a subset of tasks from fev-bench~\citep{shchur2026fevbench},  a recent forecasting benchmark that includes datasets with exogenous covariates. We restrict ourselves to those in which the target is plausibly driven by the trajectory of at least one continuous covariate, as this is the setting for which the signature method is designed. This yields four electricity load forecasting tasks (GFC12, GFC14, GFC17 and ENTSO-E), with temperature as the driving covariate, and one electricity price forecasting task (EPF), with load and generation forecasts as covariates. Each task contains several hourly series, corresponding to different zones or countries: GFC12 ($11$ series), GFC14 ($1$), GFC17 ($8$), ENTSO-E ($6$) and EPF ($5$), for a total of 31 series. 

\paragraph{Setup.} As in the previous experiment, we compare SigAdapt with Chronos-2 and TiDE. We omit GAMs as they require hand-tuning on covariates. As our setting is univariate forecasting, each series is forecasted separately with a rolling-origin protocol, two days ahead for load, to match the previous experiment, and one day ahead for prices. 
All models receive the past values of the target over the look-back window. For SigAdapt, the past signature is computed on the path of the target and the dynamic covariates. Unlike the previous setting, we do not provide any model with specific calendar variables such as holidays and days off, but we do include the classical day of the week and position in the year encodings for TiDE and SigAdapt. Similarly to the French load application, we fit one SigAdapt model per hour of the day. We do not run the TiDE model for several seeds due to computational constraints. Full details, including covariate descriptions, are given in Appendix~\ref{app:experiments_benchmarks}.

\begin{figure}[t]
\centering
\includegraphics[width=\linewidth]{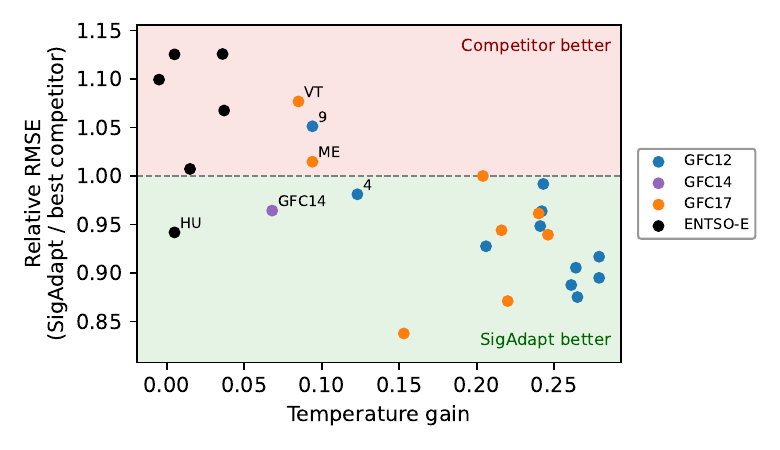}
\caption{SigAdapt's RMSE relative to the best competitor against the temperature gain, i.e.\ the increase in held-out $R^2$ obtained by adding temperature to a calendar-only linear model of the load. One point per load series.}
\label{fig:gain_vs_rmse}
\end{figure}

\paragraph{Load forecasting.}
The results are summarized in Table~\ref{tab:results_stand}, reported on standardized data. The results on the original-scale data are reported in Appendix \ref{app:experiments_benchmarks}, Table ~\ref{tab:results_raw}. SigAdapt is the best-performing model on the three GFC datasets, and it achieves the best or tied-best RMSE in 17 of the 20 series, with reductions of up to 16\% relative to the best baseline. However, it is outperformed by Chronos-2 and TiDE on the ENTSO-E dataset.  To better understand this difference in SigAdapt's performance across datasets, we use a proxy for measuring how informative temperature is in each dataset through the \emph{temperature gain}, taken as the increase in held-out $R^2$ when the current temperature and its square are added to a linear model on calendar effects to model load (Appendix~\ref{app:results_load}). Figure~\ref{fig:gain_vs_rmse} shows that SigAdapt's performance is in line with the temperature gain. It achieves the best or tied-best RMSE on every series with a gain above 0.1, and is outperformed only where temperature adds little: GFC12 zone 9, GFC17 zones ME and VT, and the ENTSO-E countries. The latter is consistent with the fact that the ENTSO-E countries rely mostly on energy sources other than electricity for space heating, as electric heating and cooling are the the main drivers of temperature sensitivity in electricity load. \footnote{In 2024, electricity accounted for a minor share of household space heating in all countries of the dataset; see Eurostat, \emph{Energy consumption in households}, \url{https://ec.europa.eu/eurostat/statistics-explained/index.php?title=Energy_consumption_in_households}.}

\begin{table}[ht!]
\centering
\caption{Test performance in terms of RMSE and MAE on standardized data. For multi-series datasets, an average is reported. Best model per row in bold.}
\label{tab:results_stand}
\resizebox{0.5\textwidth}{!}{%
\begin{tabular}{l cc | cc | cc}
\toprule
Models & \multicolumn{2}{c|}{Chronos-2} & \multicolumn{2}{c|}{TiDE} & \multicolumn{2}{c}{SigAdapt} \\
\cmidrule(lr){2-3}\cmidrule(lr){4-5}\cmidrule(lr){6-7}
Metric & RMSE & MAE & RMSE & MAE & RMSE & MAE \\
\midrule
GFC12   & 0.502 & 0.361 & 0.549 & 0.398 & \textbf{0.477} & \textbf{0.344} \\
GFC14   & 0.254 & 0.171 & 0.242 & 0.168 & \textbf{0.235} & \textbf{0.162} \\
GFC17   & 0.332 & 0.234 & 0.331 & 0.239 & \textbf{0.310} & \textbf{0.223} \\
ENTSO-E & \textbf{0.403} & \textbf{0.262} & \textbf{0.403} & 0.271 & 0.417 & 0.295 \\
EPF   & \textbf{0.393} & \textbf{0.184} & 0.460 & 0.231 & 0.399 & 0.193 \\
\bottomrule
\end{tabular}%
}
\end{table}

\paragraph{Price forecasting.}
On EPF, SigAdapt outperforms TiDE but is outperformed by Chronos-2 in all markets except PJM. Its errors are concentrated around price spikes, which a linear model on signature features captures less well than load dynamics. Improving the treatment of such jumps is a direction for future work.

\section{CONCLUSION}

We addressed the problem of forecasting a target whose dependence on the dynamics of the covariate series evolves over time, in an online setting where the model is updated as new data arrive. We proposed SigAdapt, which encodes the covariate trajectories through sliding-window signatures. The universality of these features justifies combining them with a linear model, whose dynamic effects are tracked by a Kalman filter. Since the optimal window length is unknown and may change over time, models built on different window lengths are combined by online aggregation of experts. Our experiments show that SigAdapt is most beneficial when the covariates are informative for the target, and when their trajectory matters beyond their current value, as in electricity load forecasting. We also illustrate that sliding-window signatures remain effective even in a non-stationary settings, and complement recent neural forecasters, which are difficult to update sequentially.



We note that SigAdapt is best suited to a moderate number of dynamic covariates; indeed, the signature dimension grows with the number of covariates and the truncation order. This could be addressed by dimension reduction before computing signatures. 
Also, we observed that a spike affects our signature features for as long as it remains in the window. This effect is partly mitigated by the aggregation step, since spikes leave short windows quickly. A natural direction for further robustness includes clipping or downweighting large path increments. Finally, a promising extension is multivariate forecasting, where each target series becomes a component of the signature paths used to forecast the others.


\subsubsection*{Acknowledgements}

N.D. would like to thank Adeline Fermanian for important insights and constructive feedback that helped improve this work.

\bibliography{references}

\clearpage
\appendix
\thispagestyle{empty}

\onecolumn
\aistatstitle{
Supplementary Materials}

\section{Notation glossary}
\label{app:glossary}
 
Table~\ref{tab:notation} summarizes the notation of Section~\ref{sec:method}. A subscript $t \mid T$ denotes a quantity relative to the target time $t$ that is computed from the information available at the forecast origin $T$.
 
\begin{table}[H]
\centering
\caption{Notation.}\label{tab:notation}
\small
\begin{tabular}{@{}p{0.25\linewidth}p{0.7\linewidth}@{}}
\toprule
\textbf{Symbol} & \textbf{Description} \\
\midrule
\multicolumn{2}{@{}l}{\emph{Forecasting setting}} \\
$Y_s \in \R$ & Target at time $s$ \\
$X_s \in \R^{d_X}$ & Past covariates, observed up to the forecast origin \\
$Z_s \in \R^{d_Z}$ & Future known covariates, known in advance \\
$X^{\dyn}_s$, $Z^{\dyn}_s$ & Dynamic covariates, entering a signature path \\
$X^{\point}_s$, $Z^{\point}_s$ & Pointwise covariates, entering the model through their value at a single time \\
$T$ & Forecast origin \\
$h \geq 1$ & Lead time \\
$t = T + h$ & Target (forecast) time \\
$\Info_T$ & Information available at the forecast origin $T$ \\
$\hat Y_{t\mid T}$ & Forecast of $Y_t$ issued at the origin $T$ \\
\midrule
\multicolumn{2}{@{}l}{\emph{Signatures}} \\
$\Sig^{\le N}(\cdot)$ & Signature truncated at order $N$ \\
$d$, $N$ & Dimension of a generic path and truncation order \\
$H_{[T - w,\, T]}$ & Time-augmented history path $(s, Y_s, X^{\dyn}_s, Z^{\dyn}_s)$ on $[T - w, T]$ \\
$O_{[t - w,\, t]}$ & Time-augmented outlook path $(s, Z^{\dyn}_s)$ on $[t - w, t]$ \\
$w$ & Length of the history and outlook windows \\
$N_{\mathrm H}$, $N_{\mathrm O}$ & Truncation orders of the history and outlook signatures \\
$\SH{T}$ & History signature $\Sig^{\le N_{\mathrm H}}(H_{[T - w,\, T]})$ \\
$\SO{t\mid T}$ & outlook signature $\Sig^{\le N_{\mathrm O}}(O_{[t - w,\, t]})$ \\
$\SH{T}^{(k)}$, $\SO{t\mid T}^{(\ell)}$ & Standardized level-$k$ and level-$\ell$ coordinates retained in the ridge regression \\
\midrule
\multicolumn{2}{@{}l}{\emph{Linear model}} \\
$U_{t\mid T}$ & Pointwise covariates $(X^{\point}_T, Z^{\point}_t)$ \\
$\alpha$ & Intercept \\
$\gamma$ & Coefficients of the pointwise covariates \\
$\beta^{\mathrm H}$, $\beta^{\mathrm O}$ & Coefficients of the history and outlook signatures \\
$\psi$, $\hat\psi$ & Parameters $(\alpha, \gamma, \beta^{\mathrm H}, \beta^{\mathrm O})$ and their ridge estimate \\
$\varphi_{t\mid T}$ & Features $(1, U_{t\mid T}, \SH{T}, \SO{t\mid T})$ of the linear model \\
$\lambda$ & Ridge regularization parameter \\
$\mathcal{T}_{\mathrm{tr}}$ & Set of training forecast origins \\
\midrule
\multicolumn{2}{@{}l}{\emph{Kalman filter}} \\
$\hbetaH{k}$, $\hbetaO{\ell}$ & Ridge coefficients of level $k$ of the history signature and level $\ell$ of the outlook signature \\
$\fH{k}{T}$, $\fO{\ell}{t\mid T}$ & Per-level effects: contributions of level $k$ of the history signature and level $\ell$ of the outlook signature to the batch forecast \\
$U^{\mathrm{fix}}_{t\mid T}$, $\hat\gamma_{\mathrm{fix}}$ & Pointwise covariates with coefficients kept at batch estimate, and these coefficients \\
$U^{\mathrm{ad}}_{t\mid T} \in \R^{p}$, $\hat\gamma_{\mathrm{ad}}$ & Pointwise covariates whose coefficients are adapted, and their batch estimate \\
$\xi_{t\mid T}$ & Regressor of the state-space model \\
$\theta_t$ & State at the target time $t$ \\
$\mathcal{N}(\mu, \Sigma)$ & Normal distribution with mean $\mu$ and covariance $\Sigma$ \\
$\eta_t \sim \mathcal{N}(0, Q)$ & State noise, with diagonal state noise covariance $Q$ \\
$\varepsilon_t \sim \mathcal{N}(0, \sigma^2)$ & Observation noise, with observation noise variance $\sigma^2$ \\
$\hat\theta_1$, $P_1$ & Mean and covariance of the initial state $\theta_1$ \\
$\hat\theta_{T\mid T}$ & Kalman filter estimate of the state given information up to time $T$ \\
\bottomrule
\end{tabular}
\end{table}

\section{Signatures - theory}
\label{app:sig_theory}

\subsection{Tensor spaces}
We denote by $\otimes$ the tensor product and by $(\mathbb{R}^d)^{\otimes k}$ the $k$-th tensor power of $\mathbb{R}^d$, with $(\mathbb{R}^d)^{\otimes 0} := \mathbb{R}$. Let $e_1, \ldots, e_d$ be the canonical basis of $\mathbb{R}^d$; then $\left(e_{i_1} \otimes \cdots \otimes e_{i_k}\right)_{1 \leq i_1, \ldots, i_k \leq d}$ is an orthonormal basis of $\left(\mathbb{R}^d\right)^{\otimes k}$, meaning that any element $a \in\left(\mathbb{R}^d\right)^{\otimes k}$ can be written as $a=\sum_{1 \leq i_1, \ldots, i_k \leq d} a^{\left(i_1, \ldots, i_k\right)} e_{i_1} \otimes \cdots \otimes e_{i_k}$, where $a^{\left(i_1, \ldots, i_k\right)} \in \mathbb{R}$. Furthermore, $\left(\mathbb{R}^d\right)^{\otimes k}$ is a Hilbert space of dimension $d^k$, with the following scalar product and norm:
\begin{equation*}
\langle a, b\rangle_{\left(\mathbb{R}^d\right)^{\otimes k}}=\sum_{1 \leq i_1, \ldots, i_k \leq d} a^{\left(i_1, \ldots, i_k\right)} b^{\left(i_1, \ldots, i_k\right)}, \quad \|a\|_{\left(\mathbb{R}^d\right)^{\otimes k}} =  \sqrt{\sum_{1 \leq i_1, \ldots, i_k \leq d} {\left( a^{\left(i_1, \ldots, i_k\right)} \right)} ^2} .
\end{equation*}

\begin{definition}
    We define the extended tensor algebra $T\left(\left(\mathbb{R}^d\right)\right)$ as
    \begin{equation*}
    T\left(\left(\mathbb{R}^d\right)\right):=\left\{\mathbf{a} = \left(a_0, a_1, \ldots, a_k,\ldots \right) \mid a_k \in\left(\mathbb{R}^d\right)^{\otimes k}, k \geq 0\right\}.
    \end{equation*}
\end{definition}

In other words, the extended tensor algebra is a set of infinite sequences of tensors of increasing order.
Equipped with componentwise addition and the product
\begin{equation}\label{eq:tensor_product}
    (\mathbf{a} \otimes \mathbf{b})_k := \sum_{m=0}^{k} a_m \otimes b_{k-m}, \qquad k \geq 0,
\end{equation}
$T\left(\left(\mathbb{R}^d\right)\right)$ is a non-commutative algebra with neutral element $\mathbf{1}=(1,0,0,\ldots)$. Furthermore, the subset $T_1\left(\left(\mathbb{R}^d\right)\right):=\left\{\mathbf{a} \in T\left(\left(\mathbb{R}^d\right)\right) \mid a_0=1\right\}$ is a group under $\otimes$. The following subspace will be of special interest.

\begin{definition}
    We denote with $\mathscr{T}\left(\mathbb{R}^d\right)$ the space of square-summable elements of $T\left(\left(\mathbb{R}^d\right)\right)$:
    \begin{equation*}
        \mathscr{T}\left(\mathbb{R}^d\right)=\left\{\mathbf{a} \in T\left(\left(\mathbb{R}^d\right)\right) \mid  \sum_{k=0}^{\infty}\left\|a_k\right\|_{\left(\mathbb{R}^d\right)^{\otimes k}}^2<\infty\right\}.
    \end{equation*}
\end{definition}

We can endow this space with the following scalar product and associated norm: for any $\mathbf{a}, \mathbf{b} \in \mathscr{T}\left(\mathbb{R}^d\right)$,
\begin{equation}
    \langle \mathbf{a}, \mathbf{b}\rangle_{\mathscr{T}\left(\mathbb{R}^d\right)}=\sum_{k=0}^{\infty} \langle a_k, b_k\rangle_{\left(\mathbb{R}^d\right)^{\otimes k}}, \quad
    \|\mathbf{a}\|_{\mathscr{T}\left(\mathbb{R}^d\right)} =  \sqrt{\sum_{k=0}^{\infty}\left\|a_k\right\|_{\left(\mathbb{R}^d\right)^{\otimes k}}^2} .
\label{normhilbert}
\end{equation}

The space 
$\left(\mathscr{T}\left(\mathbb{R}^d\right),\langle\cdot, \cdot\rangle_{\mathscr{T}\left(\mathbb{R}^d\right)}\right)$ is a separable Hilbert space.

\subsection{Definition and properties}\label{app:sig_definitions}

Throughout, a path is a continuous map $x : [s,t] \to \mathbb{R}^d$ of bounded variation, so that the iterated integrals below are well-defined Riemann--Stieltjes integrals.

\begin{definition}[Signature]\label{def:signature2}
    The signature of a path $x : [s,t] \to \mathbb{R}^d$ of bounded variation is the sequence of iterated integrals
    \begin{equation*}
        S\left( x_{[s,t]} \right) := \left(1, \ S^{(1)}\left( x_{[s,t]} \right), \ \ldots, \ S^{(k)}\left( x_{[s,t]} \right), \ \ldots \right) \in T_1\left(\left(\mathbb{R}^d\right)\right),
    \end{equation*}
    where, for $k \geq 1$,
    \begin{equation*}
        S^{(k)}\left( x_{[s,t]} \right) := \idotsint \displaylimits_{s < u_1<\cdots<u_k < t}  \mathrm{d} x_{u_1} \otimes \cdots \otimes \mathrm{d} x_{u_k} \in \left(\mathbb{R}^d\right)^{\otimes k}.
    \end{equation*}
    In the orthonormal basis $\left(e_{i_1} \otimes \cdots \otimes e_{i_k}\right)$ of $(\mathbb{R}^d)^{\otimes k}$, the $k$-th level $S^{(k)}\left( x_{[s,t]} \right)$ has coordinates
    \begin{equation*}
        S^{(i_1, \ldots, i_k)}\left( x_{[s,t]} \right) := \idotsint \displaylimits_{s < u_1<\cdots<u_k < t} \mathrm{d}x^{(i_1)}_{u_1} \cdots \mathrm{d}x^{(i_k)}_{u_k},
    \end{equation*}
    for $(i_1, \ldots, i_k) \in \{1, \ldots, d\}^k$, called signature coefficients. The truncated signature of order $N \geq 1$ is
    \begin{equation*}
        S^{\leq N}\left( x_{[s,t]} \right) := \left(1, \ S^{(1)}\left( x_{[s,t]} \right), \ \ldots, \ S^{(N)}\left( x_{[s,t]} \right)\right).
    \end{equation*}
\end{definition}

In applications, the truncated signature is flattened into the vector of its $s_d(N) = \sum_{k=0}^{N} d^k$ coefficients, ordered by level and lexicographically within each level.

\begin{proposition}[Properties of the signature]\label{prop:sig_properties}
    Let $x : [s,t] \to \mathbb{R}^d$ and $y : [t,u] \to \mathbb{R}^d$ be paths of bounded variation.
    \begin{enumerate}
        \item (\textbf{Invariance to reparametrization}) For any continuously differentiable non-decreasing surjection $\psi : [s,t] \to [s,t]$, the path $\widetilde{x}_v := x_{\psi(v)}$ satisfies $S\left( \widetilde{x}_{[s,t]} \right) = S\left( x_{[s,t]} \right)$.
        \item (\textbf{Invariance to translation}) For any $a \in \mathbb{R}^d$, the path $\overline{x}_v := x_v + a$ satisfies $S\left( \overline{x}_{[s,t]} \right) = S\left( x_{[s,t]} \right)$.
        \item (\textbf{Chen's identity}) Let $x * y : [s,u] \to \mathbb{R}^d$ denote the concatenation of $x$ and $y$, defined by $(x*y)_v := x_v$ for $v \in [s,t]$ and by $(x*y)_v := x_t + y_v - y_t$ for $v \in [t,u]$. Then
        \begin{equation*}
            S\left( (x * y)_{[s,u]} \right) = S\left( x_{[s,t]} \right) \otimes S\left( y_{[t,u]} \right).
        \end{equation*}
        \item (\textbf{Time reversal}) The time-reversed path $\overleftarrow{x}_v := x_{s+t-v}$ satisfies
        \begin{equation*}
            S\left( x_{[s,t]} \right) \otimes S\left( \overleftarrow{x}_{[s,t]} \right) = \mathbf{1},
        \end{equation*}
        that is, $S\left( \overleftarrow{x}_{[s,t]} \right)$ is the inverse of $S\left( x_{[s,t]} \right)$ in the group $T_1\left(\left(\mathbb{R}^d\right)\right)$.
        \item (\textbf{Factorial decay}, \citealp{Lyons_2007}) For all $k \geq 0$,
        \begin{equation*}
            \big\| S^k\left( x_{[s,t]} \right) \big\|_{(\mathbb{R}^d)^{\otimes k}} \leq \frac{\|x\|_{1\text{-var},[s,t]}^k}{k!}.
        \end{equation*}
        In particular, $S\left( x_{[s,t]} \right) \in \mathscr{T}(\mathbb{R}^d)$.
    \end{enumerate}
\end{proposition}

Properties 1 and 2 show that the signature does not encode the speed at which a path is traversed and its starting point. Adding time as a coordinate, $u \mapsto (u, x_u)$, fixes the invariance to reparametrization. Prepending the origin as a basepoint removes the remaining translation invariance.

\begin{example}[Linear path]\label{ex:piecewise_2}
    Let $x : [s,t] \to \mathbb{R}^d$ be the linear path $u \mapsto x_s + \frac{u-s}{t-s}\,(x_t - x_s)$. Then, for all $k \geq 0$,
    \begin{equation*}
        S^k\left( x_{[s,t]} \right) = \frac{1}{k!}\,(x_t - x_s)^{\otimes k}.
    \end{equation*}
    Combined with Chen's identity, the signature of a piecewise linear path is therefore the tensor product of the signatures of its linear segments, which is how it is computed in practice, e.g. in the Python package \texttt{iisignature}~\citep{iisignature}.
\end{example}

\section{Signatures - fast updates}
\label{app:sig_algo}

To compute signatures efficiently online, we exploit the large overlap between consecutive windows. Instead of recomputing the signature at each step, we update it using the properties of Section~\ref{app:sig_theory}. This extends the sliding-window update of \citet{drobac2025sliding} to basepoint-augmented paths, i.e.\ paths that always start at zero. 

When the window moves from $[t-w-1, t-1]$ to $[t-w, t]$, the update consists of three tensor products:
\begin{enumerate}
    \item \textbf{Remove the start of the window.} The first two segments of the path are the basepoint segment and the oldest segment. By time reversal, the signature of these segments traversed backwards is their inverse. Multiplying it on the left removes them.
    \item \textbf{Add the new segment.} By Chen's identity, we multiply the current signature on the right by the signature of the new segment.
    \item \textbf{Add the new basepoint.} By Chen's identity, we multiply the current signature on the left by the signature of the segment from zero to the new first point.
\end{enumerate}

The summary in pseudocode can be found in Algorithm \ref{algo:update_sig}).

For numerical stability, we standardize the covariates, so that the basepoint segment is not too long. 

Finally, we note that signatures are computed with the Python package \texttt{iisignature}~\citep{iisignature}. The sliding-window updates chain the functions \texttt{sigcombine} and \texttt{sigjoin}, whose outputs are returned in single precision in the released version. We therefore use the \texttt{doubles} branch of the package, which returns them in double precision\footnote{\url{https://github.com/bottler/iisignature/pull/25}}. This prevents rounding errors from accumulating over successive updates. To avoid the latter, we also recompute the signature from scratch every 500 observations.

\begin{algorithm}
\caption[Basepoint signature update on sliding windows]{Basepoint signature update on sliding windows. Here $\ell_{a,b}$ denotes the linear segment from $a$ to $b$, with $S^{k}(\ell_{a,b}) = \frac{1}{k!}(b-a)^{\otimes k}$.}
\label{algo:update_sig}
\begin{algorithmic}[1]
    \REQUIRE data $\left\{ x_t \mid t \geq 0 \right\}$, window size $w$, truncation order $N$
    \ENSURE $S^{\leq N}\!\left(\widehat{x}_{[t-w,t]}\right)$ for all $t \geq w$
    \STATE Compute initial signature $S \gets S^{\leq N}\!\left(x_{[0,w]}\right)$
    \STATE Add the basepoint: $S^{\leq N}\!\left(\widehat{x}_{[0,w]}\right) \gets S^{\leq N}\!\left(\ell_{0,x_0}\right) \otimes S$
    \FOR{$t = w+1, w+2, \ldots$}
        \STATE $S_{\text{old}} \gets S^{\leq N}\!\left(\overleftarrow{x}_{[t-w-1,t-w]}\right)$
            \COMMENT{reversed segment}
        \STATE $S \gets S_{\text{old}} \otimes S$
            \COMMENT{time reversal}
        \STATE $S_{\text{new}} \gets S^{\leq N}\!\left(x_{[t-1,t]}\right)$
            \COMMENT{new segment}
        \STATE $S \gets S \otimes S_{\text{new}}$
            \COMMENT{Chen's identity}
        \STATE $S^{\leq N}\!\left(\widehat{x}_{[t-w,t]}\right) \gets S^{\leq N}\!\left(\ell_{0,x_{t-w}}\right) \otimes S$
            \COMMENT{basepoint}
    \ENDFOR
\end{algorithmic}
\end{algorithm}

\section{Experiments - Case study}
\label{app:french_load}

\subsection{Setup}

\paragraph{Data.} The dataset contains the French national electricity consumption at a half-hourly resolution, published by RTE, together with a national temperature series and calendar information (public holidays and school holidays for the three French school zones A, B and C), from 1 January 2018 to 30 June 2026. 
The period 2018 - 2023 is used for training, and 1 January 2024 to 30 June 2026 for testing ($43\,776$  data points). Hyperparameters are selected on 2022 - 2023 with models fitted on 2018 - 2021, and all models are then refitted on 2018 - 2023 before being evaluated on the test period. The dataset, including the covariates, is be made public with the paper.

\paragraph{Forecasting task.} The target $Y_t$ is the national consumption, forecast with a lead time of $48$ hours (2 days). At forecast origin $T$, to forecast $Y_{T+\text{2 days}}$, the input of the models is the past consumption, up to time $T$,  while temperature and calendar variables are treated as future covariates, available up to $T+\text{ 2 days}$. Realised temperatures are used in place of temperature forecasts, for all models. 

\paragraph{Calendar covariates.} All models share the calendar variables in Table~\ref{tab:rte_calendar}. A day is a working day if it is a Monday to Friday and not a public holiday. A bridge is a working day preceded and followed by a non-working day; bridges falling within the end-of-year break are not flagged, so that their effect does not add to the Christmas effects. The Christmas indicator covers 24 December to 1 January, on weekdays only, again to avoid the addition of the weekend effect and the holiday effect. School holidays are considered only when they hold in all three zones, as we are forecasting on a national level. The Christmas school holidays are then restricted to the days outside 24 December to 1 January, so that the two indicators are disjoint, and the summer school holidays are split into July and August, as the latter also includes the more general summer shutdown in France.

\begin{table}[h]
\centering
\caption{Calendar covariates for RTE. All indicators are evaluated per day and repeated over its $48$ half-hours.}
\label{tab:rte_calendar}
\begin{tabular}{ll}
\toprule
Variable & Definition \\
\midrule
Position in the year & $\sin$ and $\cos$ of $2\pi \times$ (fraction of the year elapsed) \\
Day of the week & dummies, Monday as reference \\
Working day & Monday--Friday and not a public holiday \\
Bridge & working day between two non-working days, outside the Christmas period \\
Christmas & 24~Dec.--1~Jan., Monday--Friday \\
Christmas school holidays & school holidays (all zones) in Dec.--Jan., outside 24~Dec.--1~Jan. \\
July holidays & school holidays (all zones) in July \\
August holidays & school holidays (all zones) in August \\
\bottomrule
\end{tabular}
\end{table}

\paragraph{Evaluation metrics.}
All reported results are computed on the test period. We evaluate point forecasts using the root mean squared error (RMSE), the mean absolute error (MAE) and the mean absolute percentage error (MAPE). Denoting by $T$ the number of data points in the test period, by $Y_t$ the observed value and by $\hat{Y}_t$ its forecast, these are defined as 
\begin{equation}
    \mathrm{RMSE} = \sqrt{\frac{1}{T}\sum_{t=1}^T \big( \hat{Y}_t-Y_t \big)^2}, \quad \mathrm{MAE} = \frac{1}{T}\sum_{t=1}^T  |\hat{Y}_t-Y_t|, \quad \mathrm{MAPE} = \frac{100}{T}\sum_{t=1}^T\frac{|\hat{Y}_t-Y_t|}{|Y_t|}.
    \label{eq:metrics}
\end{equation}

Since the MAPE becomes unstable when some observed values are close to zero, we report it only for datasets where this does not occur, and omit it otherwise (EPF, GFC12-4, GFC12-9).

\subsection{Models}
\label{app:french_load_models}

\paragraph{SigAdapt configuration.} A separate model is fitted for each of the $48$ half-hours of the day. For a window of length $w$, the model combines two signatures truncated at depth $6$: the history signature of the path consisting of time, temperature and load, whose window ends at the forecast origin, and the outlook signature of the time-augmented temperature path, whose window ends at the forecast time. Signature components involving only time, or only time and temperature, are removed from the history signature, and components involving only time are removed from the outlook signature. The pure-time components are constant and would only add an intercept to the regression model, while the time-and-temperature components are already contained in the outlook signature.
Temperature and consumption are standardized on the training period, and the signature terms are standardized for each half-hour before entering the ridge regression. The calendar variables enter the ridge regression together with the signature terms, without standardization. The ridge regularization parameter is selected on the 2022--2023 validation period. The bridge indicator is left essentially unpenalized (its penalty is divided by $10^{4}$): with about three bridges per year, a shared penalty would shrink most of its effect. The fitted signature terms are then summed within each signature level, giving one feature per signature and level. A Kalman filter is then fitted, in which only the intercept and the level coefficients, termed signature dynamic effects, are allowed to drift. Its covariance matrices are selected on the training period by the iterative grid search of~\citet[Chapter 5]{these_joseph}.
The filter is updated with a delay of two days, matching the lead time, and is not updated on special days (public holidays falling on a weekday, bridges, and 24 December to 1 January), whose residuals are not representative of the level of normal days. This procedure is repeated for the windows $w \in \{2, 7, 14, 28, 84, 168\}$ days, and the resulting forecasts are combined for each half-hour by online expert aggregation (BOA, \texttt{opera} package). The weights are updated using data up to $t$ to forecast $t + 2$ days, again in line with the two-day-ahead setting. 

\paragraph{GAM and GAMAdapt configuration.} As a classical benchmark for load forecasting, we fit one generalized additive model per half-hour, using the Python package \texttt{pygam}. Our GAM includes: cyclic P-spline in the position in the year, splines in temperature, in two exponentially smoothed temperatures (smoothing factors $0.01$ and $0.05$), and the past load from two and seven days before, as well as factor terms for the day of the week and for each calendar indicator of Table~\ref{tab:rte_calendar}. The smoothing parameters are selected by Generalized Cross Validation (GCV) score, shared across terms, except for the bridge, whose penalty is again divided by $10^{4}$. In the adaptive version, named GAMAdapt, a Kalman filter is fit on the partial effects following~\citet{obst2021adaptive}) of the temperature, the two smoothed temperatures, lagged loads and the intercept. All other terms, including the calendar effects and the are held fixed. The filter is tuned, delayed and paused on special days exactly as for SigAdapt. We report both the static GAM and GAMAdapt.

\paragraph{TiDE configuration.} TiDE~\citep{das2024tide} is an MLP-based encoder-decoder deep neural forecaster, for which we use the implementation of the \texttt{darts} library~\citep{herzen2022darts} (version 0.28.0). As input, it receives $18$ covariates: temperature, the calendar variables of Table~\ref{tab:rte_calendar} (with all seven day-of-week dummies), and the sine and cosine of the time of day, which replaces the per-half-hour structure of SigAdapt and the GAM. The model outputs predictions over the $96$ forecasting horizon (2 days of half hourly data) and is used with a rolling origin moved by one half-hour, keeping only the last point of each output, so that every half-hour is forecast exactly $48$ hours ahead. TiDE is the only model whose training is random and depends on the seed; and as a result so do the forecasts. We therefore train it with five seeds, $\{1, 42, 145, 150, 1000\}$, and report the mean and standard deviation of the test metrics over the seeds. The look-back is selected among $\{1, 2, 4, 7, 14\}$ days on the validation set. To limit the cost of the $25$ training runs, validation forecasts use a stride of $47$ half-hours, which covers all times of the day and days of the week evenly at the same lead time. Based on the best RMSE averaged over the five seeds, the look-back window is set to 7 days. The final model is refitted on the joint train and validation period for each seed. The remaining configuration hyperparameters can be found in the code.

\paragraph{Chronos-2 configuration.} Chronos-2~\citep{ansari2025chronos2} is a pretrained foundation model that is applied zero-shot. It receives temperature and the special-day indicators (public holidays falling on a weekday, bridges, Christmas, Christmas school holidays, July and August holidays) as past and future covariates, but no day-of-week, time-of-day or seasonal encoding, so that the daily, weekly and yearly patterns are inferred from the context. Forecasts are produced with the same protocol as for TiDE: an output of $96$ half-hours, a rolling origin moved by one half-hour, and only the last point kept. The point forecast is the median of the predictive distribution. The context length is selected among $\{2, 7, 14, 28, 84, 168\}$ days on 2022 - 2023, with a stride of one half-hour. The selected context is $168$ days ($8\,064$ half-hours), which is the largest value considered on the grid, as longer contexts would approach the maximum context length of the model. This long context, however, makes the forecasting quite slow.

\subsection{Ablation study}\label{app:rte_ablation}

 \paragraph{Setup.} To assess the role of each component of SigAdapt (signature features, Kalman filter and aggregation over windows), we compare it with variants in which one or several components are removed. All other choices are kept identical: one ridge model per half-hour with the calendar variables of Table~\ref{tab:rte_calendar}, the ridge parameter selected on 2022--2023, the same Kalman tuning, and the same aggregation. We consider three scenarios:
\begin{enumerate}[label=\roman*)]
    \item \emph{Without signatures.} A linear model uses the temperature at the target time $t$, the load at $t - 2$ days and the calendar variables. It is evaluated without and with the Kalman filter, in which the intercept, the temperature and load coefficients drift. To mirror the windows of SigAdapt, we then build six such linear experts, each adding the temperature $w$ days before the target, for $w \in \{2, 7, 14, 28, 84, 168\}$, and aggregate them as in SigAdapt.
    \item \emph{Without the Kalman filter.} The six static signature experts are aggregated directly.
    \item \emph{Without aggregation.} A single signature expert with the Kalman filter is used, with its window chosen with hindsight on the test period: either one window for all half-hours, or the best window for each half-hour. Since these windows are selected on the test period itself, they bound the accuracy of any procedure that selects a single window, for instance on a validation period.
\end{enumerate}
Finally, we include GAMAdapt, which shares the calendar variables, the information set and the Kalman filter with SigAdapt, but represents temperature and load through smooth functions of the current and exponentially smoothed temperatures and of lagged loads, without aggregation.

\paragraph{Results.} Table~\ref{tab:rte_ablation_full} reports the test errors. Without signatures, the linear model reaches an RMSE of $2\,632$\,MW without the Kalman filter and $1\,625$\,MW with it, and aggregating the six linear experts barely improves on the best of them ($1\,607$ vs.\ $1\,606$\,MW). This is mainly because the errors of the six linear experts are almost perfectly correlated (mean pairwise correlation $0.987$, at least $0.965$ for every pair): each expert differs from the others by a single lagged temperature, and they make essentially the same mistakes. The errors of the signature experts are much less correlated ($0.61$ on average). Neighbouring windows are the most correlated ($0.76$ between $14$ and $28$ days, $0.75$ between $84$ and $168$ days) and distant windows the least ($0.50$ between $2$ and $28$ days), indicating that each window summarises the temperature and load paths at its own time scale. With signatures, the best static expert reaches $1\,566$\,MW; adding either the Kalman filter (best window per half-hour: $1\,257$\,MW) or the aggregation ($1\,260$\,MW) reduces the error by about $20\%$, and combining both, as in SigAdapt, gives $1\,168$\,MW. GAMAdapt ($1\,280$\,MW) is slightly more accurate than the best single signature window ($1\,288$\,MW), and less accurate than SigAdapt.

Signatures provide two things that raw inputs do not. First, nonlinear features of the temperature and load paths: a single signature window reaches the accuracy of GAMAdapt, whose temperature and load features hand-crafted, and is $21\%$ more accurate than the linear model with the same information. Part of this gain reflects the nonlinear response of load to temperature, which a model linear in temperature cannot represent. Second, diversity across windows, which is what makes aggregation effective: the same aggregation applied to the linear experts brings no gain. The aggregation  can partly compensate for the absence of the Kalman filter by combining static experts online. SigAdapt combines both forms of adaptation, and is more accurate than any variant in which one of the three components is removed.

\begin{table}[ht]
\centering
\small
\begin{tabular}{lccccc}
\toprule
Variant & Kalman & Aggr. & RMSE & MAE & MAPE \\
\midrule
\multicolumn{6}{l}{\emph{Without signatures}} \\
Linear model & no & no & 2\,632 & 2\,004 & 3.99 \\
Linear model & yes & no & 1\,625 & 1\,188 & 2.36 \\
Best linear expert$^{\dagger}$ & yes & no & 1\,606 & 1\,180 & 2.34 \\
Linear experts & yes & yes & 1\,607 & 1\,178 & 2.34 \\
\midrule
\multicolumn{6}{l}{\emph{With signatures}} \\
Best static expert$^{\dagger}$ & no & no & 1\,566 & 1\,131 & 2.20 \\
Static experts & no & yes & 1\,260 & 902 & 1.75 \\
Best single window$^{\dagger}$ & yes & no & 1\,288 & 947 & 1.88 \\
Best window per half-hour$^{\dagger}$ & yes & no & 1\,257 & 926 & 1.83 \\
SigAdapt & yes & yes & \textbf{1\,168} & \textbf{845} & \textbf{1.66} \\
\midrule
\multicolumn{6}{l}{\emph{Hand-crafted features}} \\
GAM & no & no & 1\,504 & 1\,111 & 2.16 \\
GAMAdapt & yes & no & 1\,280 & 939 & 1.85 \\
\bottomrule
\end{tabular}
\caption{Ablation on French national load, test period 1 January 2024 to 30 June 2026; RMSE and MAE in MW, MAPE in \%. $^{\dagger}$Chosen with hindsight on the test period.}
\label{tab:rte_ablation_full}
\end{table}

\subsection{Additional graphs and tables}

\begin{figure}[H]
\centering
\includegraphics[width=0.8\linewidth]{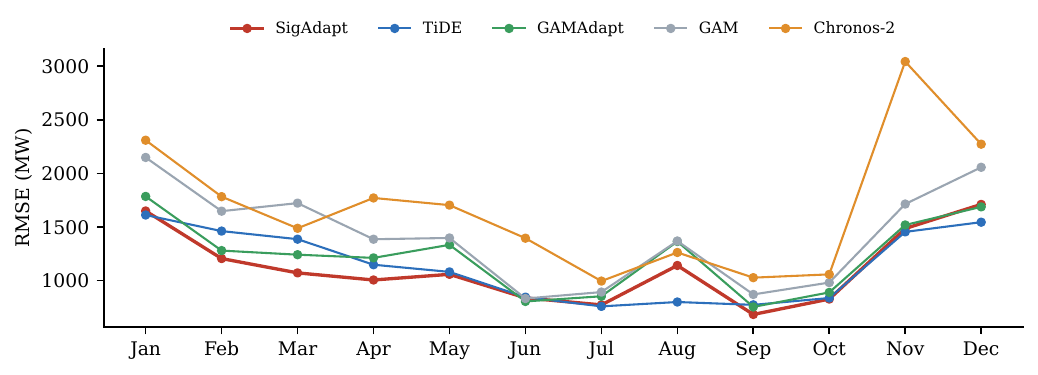}
\caption{RMSE on the test set (January 2024 - June 2026) by month for all models.}
\label{fig:rmse_by_month}
\end{figure}

Figure~\ref{fig:rmse_by_month} shows the RMSE of each model per calendar month, pooling the test years. TiDE outperforms SigAdapt mainly in August and December, which suggests that it better captures the change in load during holiday periods.

\section{Experiments - Benchmark datasets}
\label{app:experiments_benchmarks}

\subsection{Experimental setup}
\label{app:setup}

This appendix provides a detailed description of the experiments reported in Section~\ref{sec:bench_data}. For each dataset, we describe the data,
the forecasting task, and the model configuration. 
Results are reported in Section~\ref{app:results_load} for the electricity load datasets and in Section~\ref{app:results_epf} for the electricity price dataset.

\paragraph{Models.} As before, we compare SigAdapt with the two neural-network-based competitors: Chronos-2  and TiDE. We omit the comparison with the GAM as it was a specialized model, and our goal on the benchmark datasets is to apply the procedure with minimal fine-tuning.

\paragraph{Data.} We evaluate our method on a subset of tasks from fev-bench~\citep{shchur2026fevbench}. We restrict ourselves to tasks in which the target is plausibly driven by the trajectory of at least one continuous covariate, as this is the setting that motivates the use of signature features. This yields four load forecasting tasks (GFC12, GFC14, GFC17 and ENTSO-E) and one electricity price forecasting task (EPF), all at hourly frequency. All data were accessed through the \texttt{fev} Python package\footnote{\url{https://github.com/autogluon/fev}}. We use the task definitions provided by the benchmark for the target and for the known and past dynamic covariates, but depart from its evaluation protocol to match our setting and the RTE application.

\paragraph{Forecasting protocol.} All tasks are treated as univariate point forecasting problems. As most of the datasets contain several series (e.g. zones, countries or markets), each target series is forecast separately, using only its own history and its own covariates. 
We use a rolling-origin forecasting protocol in which the forecast origin advances by one time step. At each origin $T$, every method uses the information available at time $T$, together with the values of known-future covariates up to $T+h$, to forecast the target at time $T+h$, where $h$ denotes the lead time. 
The origin is then moved to $T+1$ and the procedure is repeated. Since Chronos-2 and TiDE produce multi-step forecasts, we set their forecast horizon to $h$ and retain only the final step, so that all methods are evaluated on the same target timestamps at the same lead time. For all the electricity load datasets, we consider two-day-ahead forecasts, i.e. $h=2 \text{ days}$ to match our operational setting, and for the electricity price we consider one-day-ahead forecasts to match other literature on this dataset. 
We also note that for all of the tasks including temperature as a covariate, realized future temperatures are used in place of temperature forecasts, as is common practice in this type of forecasting framework.

\paragraph{Data splits and model selection.} Each dataset is split chronologically into training, validation and test periods. The validation period is used only for model selection and hyperparameter tuning. For all models, the amount of history is taken from the same grid of $\{2, 7, 14, 28, 56, 84\}$~days. For Chronos-2 and TiDE, a single value is selected on the validation period, defining the context length and the look-back window, respectively. For SigAdapt, no single window is selected; instead, the model fitted with each window length is treated as an expert, and the experts are combined by online aggregation. The validation period is instead used to tune the ridge penalty, while the joint train and validation period is used to select the Kalman filter hyperparameters by iterative grid search, using the \texttt{viking} package\footnote{\url{https://pypi.org/project/viking-kalman/}}. After model selection, TiDE and SigAdapt are refitted on the combined training and validation periods, while Chronos-2 is applied zero-shot and requires no fitting.

\paragraph{Covariates.} For every task, the models receive the past values of the target over the look-back window. In addition, the load tasks use temperature as a covariate; for ENTSO-E, direct and diffuse solar radiation are also included. For EPF, the covariates are the load and generation forecasts. In all datasets, calendar variables carry substantial predictive information, with the hour of the day being the strongest signal. SigAdapt accounts for the hour of the day by being fitted separately for each hour, and additionally receives the sine and cosine of the position in the year and day-of-the-week dummy variables. For the two competitor models, we follow the choices made in their original papers. TiDE receives the same calendar covariates as SigAdapt, together with the sine and cosine of the time of day, so that it has access to the same information~\citep[App.~B.2]{das2024tide}. 
Chronos-2 receives no calendar covariates. For SigAdapt, Table~\ref{tab:routing} summarizes which variables enter the past signature, the outlook signature, and the linear part of the model for each task.

\paragraph{Evaluation metrics.}
All reported results are computed on the test period. We again evaluate point forecasts using RMSE, MAE and MAPE, as defined in Equation \ref{eq:metrics}. Since the MAPE becomes unstable when some observed values are close to zero, we report it only for datasets where this does not occur, and omit it otherwise (EPF, GFC12-4, GFC12-9).

\begin{table}[ht]
\centering
\caption{Covariate routing in SigAdapt for each task. The target signature is computed over the look-back window ending at the forecast origin; the covariate signature over a window ending at the forecast time. The remaining variables enter the ridge regression directly. The target is in bold.}
\label{tab:routing}
\begin{tabular}{l lll}
\toprule
Task & History signature & Outlook signature & Pointwise features \\
\midrule
GFC12, GFC14, GFC17 & time, temperature, \textbf{load} & time, temperature & calendar encodings \\
ENTSO-E & time, temperature, \textbf{load}  & time, temperature & calendar encodings, direct  \\ 
& & & and diffuse radiation \\
EPF & time, load, generation, \textbf{price}$^{*}$ & time, load, generation & calendar encodings \\
\bottomrule
\end{tabular}
\vspace{4pt}
\begin{minipage}{\linewidth}
\footnotesize
$^{*}$For EPF, the history signature is truncated at level 4, and only its components involving the price are kept, that is, the pure price terms and the cross terms between the price and the covariates, since the pure covariate terms are included in the outlook signature.
\end{minipage}
\end{table}

\subsection{Tasks and configurations per dataset}
\label{app:datasets}

\subsubsection{GFC12}

\paragraph{Data.} The dataset originates from the Global Energy Forecasting Competition 2012~\citep{GFC12_Hong_2014}, but we emphasize that we use the version provided by fev-bench. It contains hourly load and temperature series for $11$ zones, spanning January 2004 to June 2008. The data are split into training (2004--2005), validation (2006) and test (January 2007 to December 2007) periods. For zone~10, we omit the last six months of the original dataset because of a large level shift in the load. This period was not part of the test data in the original competition and has often been treated separately in the literature \citep[see e.g.][Sec. 3.6]{nedellec2014gefcom12}.

\paragraph{Task and configuration.} For each of the zones, the target is the load and the covariate is temperature. Temperature is treated as a known covariate, meaning that its values are available both over the context window and the forecast horizon. As a continuous variable that drives the load, temperature enters the SigAdapt model through the signature of its path.

\subsubsection{GFC14}

The dataset originates from the Global Energy Forecasting Competition 2014~\citep{hong2016gefcom2014}, but we again use the fev-bench version for load. It contains a single hourly load series with an associated temperature series, covering a two year period, from 01/01/2013 to 31/12/2014. The first year is used for training, the following six months for validation, and the last six months for testing. The target is the load and the covariate is temperature, which again enters the SigAdapt model through the signature.

\subsubsection{GFC17}

The dataset is the fev-bench version of the Global Energy Forecasting Competition 2017~\citep{hong2019gefcom2017} one. It contains hourly load and temperature series for $8$ zones, again spanning two years, from 01/01/2016 to 31/12/2017. The first year is used for training, the following six months for validation, and the last six months for testing.
As for GFC12 and GFC14, the target is the load and the covariate is temperature, entering SigAdapt through the signature.

\subsubsection{ENTSO-E}

\paragraph{Data.} We use the hourly version of the ENTSO-E load dataset from fev-bench, which contains load series for $6$ countries from 2015 to 2019. The first three years are used for training, 2018 for validation, and 2019 for testing.

\paragraph{Task and configuration.} The target is the national load, and the covariates are temperature, diffuse horizontal radiation and direct horizontal radiation, all treated as known-future covariates. In SigAdapt, temperature enters through the signature path, as its recent dynamics influence the load, whereas the two radiation covariates enter the ridge regression directly. Chronos-2 and TiDE receive all three covariates both over the context and as known-future covariates.

\subsubsection{EPF}

\paragraph{Data.} The EPF dataset~\citep{lago2021forecasting} contains hourly day-ahead electricity prices for five markets: the French (FR), German (DE) and Belgian (BE) EPEX markets, Nord Pool (NP), and the Pennsylvania--New Jersey--Maryland (PJM) market. Each series contains six years of data, with the exact spans given in Table~\ref{tab:epf_data}. The first three years are used for training, the fourth for validation, and the last two for testing, following the recommendation of~\citet{lago2021forecasting}. Unlike load, prices exhibit sharp spikes and can take values close to or below zero, which motivates the transformation described below and the exclusion of the MAPE for this dataset.

\paragraph{Forecasting task.} The target is the hourly day-ahead price. Each market comes with two known covariates, listed in Table~\ref{tab:epf_data}, which are day-ahead forecasts published by the system operators rather than realised values.

\paragraph{Model configuration.} For SigAdapt, we apply a standard variance-stabilising transformation for electricity price data~\citep{uniejewski2018variance}: the target and the covariates are standardised by their median and median absolute deviation and then transformed with the inverse hyperbolic sine (asinh), as implemented in the \texttt{epftoolbox} library~\citep{lago2021forecasting}. For the Belgian market, only the median standardisation is applied to the target. Forecasts are mapped back to the original scale for evaluation. Both covariates enter SigAdapt through the signature, over both the past and the future windows. Chronos-2 receives untransformed data, since it applies an asinh-based transformation internally to each input series~\citep[Sec.~3.1]{ansari2025chronos2}. For TiDE, reversible instance normalisation~\citep{kim2021reversible} is enabled.

\begin{table}[ht]
\centering
\caption{Dataset spans and known covariates for EPF markets.}
\label{tab:epf_data}
\begin{tabular}{llll}
\toprule
Market & Period & Covariate 1 & Covariate 2 \\
\midrule
FR  & 09/01/2011 -- 31/12/2016 & Generation forecast & System load forecast \\
DE  & 09/01/2012 -- 31/12/2017 & Load forecast (Amprion zone) & Wind and solar generation forecast \\
BE  & 09/01/2011 -- 31/12/2016 & Generation forecast (France) & System load forecast (France) \\
NP  & 01/01/2013 -- 24/12/2018 & Grid load forecast & Wind generation forecast \\
PJM & 01/01/2013 -- 24/12/2018 & System load forecast & Load forecast (COMED zone) \\
\bottomrule
\end{tabular}
\end{table}

\subsection{Results on the load datasets}
\label{app:results_load}

\begin{figure}[H]
    \centering
\includegraphics[width=0.9\linewidth]{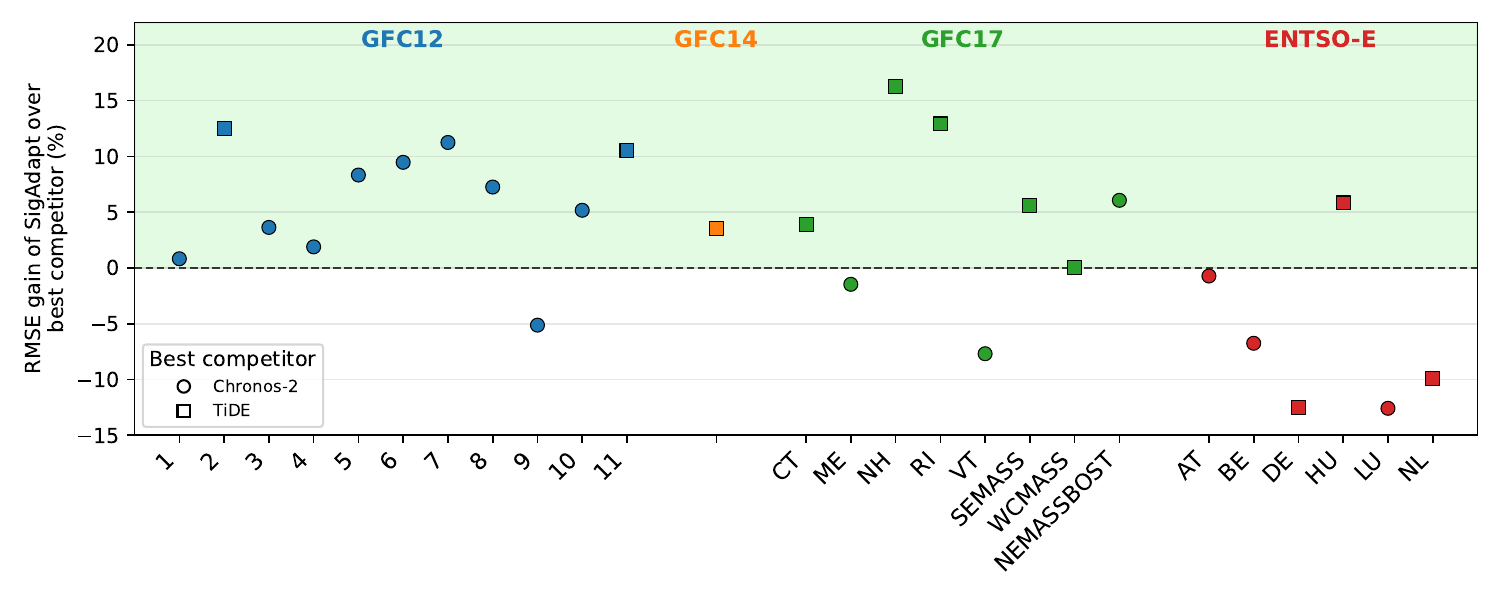}
    \caption{Relative RMSE gain of SigAdapt over the best competitor for each
forecasted series, computed as
$100 \times (\mathrm{RMSE}_{\mathrm{best}} - \mathrm{RMSE}_{\mathrm{SigAdapt}}) / \mathrm{RMSE}_{\mathrm{best}}$,
where $\mathrm{RMSE}_{\mathrm{best}} = \min(\mathrm{RMSE}_{\text{Chronos-2}}, \mathrm{RMSE}_{\text{TiDE}})$.
Series are grouped by dataset along the horizontal axis. Marker shape
indicates which competitor attains the lower RMSE. Positive values (shaded
region) indicate that SigAdapt outperforms both competitors.}
    \label{fig:rmse_gain_load}
\end{figure}

Results on the original scale of the data are reported per series in
Table~\ref{tab:results_raw} and summarized in Figure~\ref{fig:rmse_gain_load}
by the relative RMSE gain of SigAdapt over the best competitor.

On GFC12, SigAdapt outperforms both competitors on all zones and for all metrics, except for Zone~9, where it is outperformed by Chronos-2, and for the MAPE in Zone~1. 
On GFC14, SigAdapt outperforms both Chronos-2 and TiDE on all metrics.
On GFC17, SigAdapt achieves the lowest RMSE in five of the eight zones and ties with TiDE in WCMASS, with reductions of up to 16\% relative to the best baseline (NH). Since the differences with Chronos-2 in ME are small, VT is the only zone with a substantial gap.
On ENTSO-E, SigAdapt achieves the lowest RMSE only in HU. 

\begin{figure}[ht]
    \centering
    \includegraphics[width=0.9\linewidth]{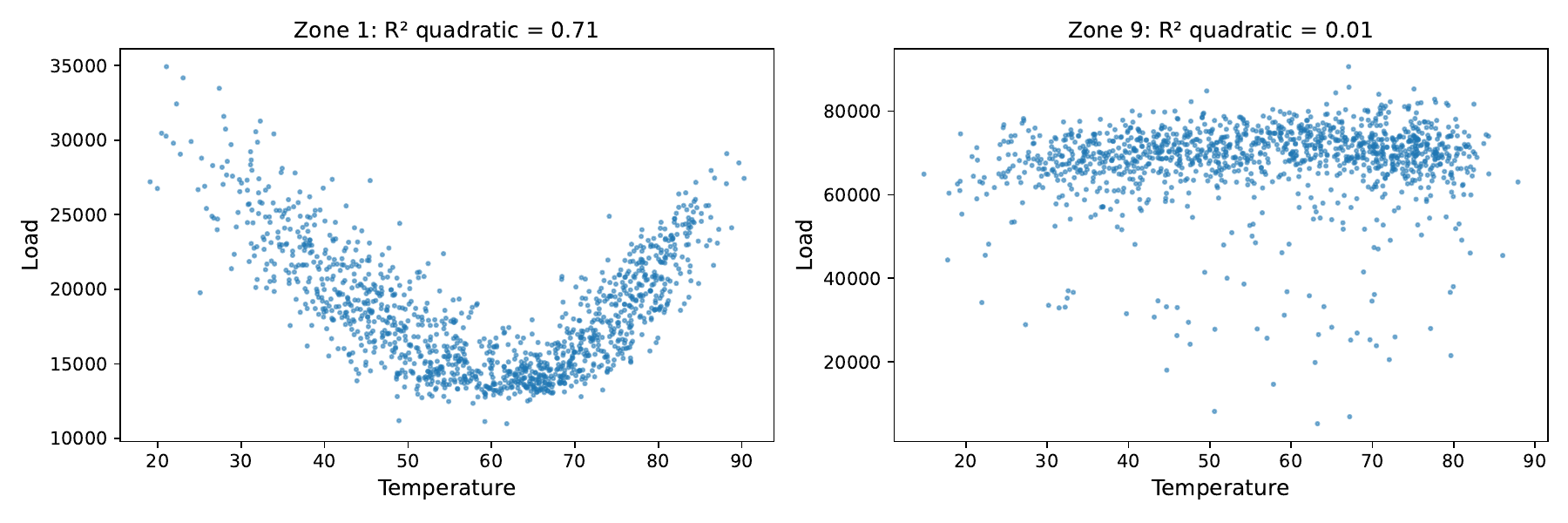}
    \caption{Daily mean load against daily mean temperature over the training period, for Zone 1 (left) and Zone 9 (right). Zone 1 shows the typical quadratic relationship, with load higher on both cold and hot days; a quadratic regression of load on temperature gives $R^2 = 0.71$. In Zone 9, load is essentially insensitive to temperature, with $R^2 = 0.01$.}
    \label{fig:temp_zones19}
\end{figure}

This variation can in part be explained by the fact that SigAdapt is designed to exploit a continuous exogenous signal that drives the target, by capturing this influence through the signature. When the covariate carries no information about the target, the signature terms involving it add variance without adding predictive power. 
This is precisely the situation in Zone~9 of GFC12, whose load has been reported to show no clear dependence on temperature \citep[see e.g.][Sec.~3.6]{nedellec2014gefcom12}. Figure~\ref{fig:temp_zones19} illustrates this: a quadratic regression of daily load on daily temperature explains most of the variation in Zone~1 ($R^2 = 0.71$) but almost none in Zone~9 ($R^2 = 0.01$).

\paragraph{Informativeness of temperature.}
To understand the variation in SigAdapt's performance across series, we use linear models as a simple but principled proxy for how informative temperature is for load prediction, with modeling choices that follow standard practice in load forecasting. The relationship between load and temperature is nonlinear and U-shaped, reflecting heating demand at low temperatures and cooling demand at high temperatures \citep{bessec2008nonlinear}, which we capture by a quadratic model. Moreover, load responds to temperature with a delay, owing for instance to the thermal inertia of buildings. This effect is often captured through exponentially smoothed temperatures in additive load models \citep{nedellec2014gefcom12, goude2013local, gaillard2016additive}.

For each series, we fit three nested linear regressions of the hourly load, by least squares, on: (i) calendar effects, namely hour-of-day and day-of-week dummy variables, as well as the sine and cosine of the position in the year; (ii) the current temperature $T_t$ and its square; and (iii) the exponentially smoothed temperatures $S^\alpha_t = \alpha T_t + (1-\alpha)S^\alpha_{t-1}$ and their squares, for $\alpha \in \{0.05, 0.01\}$. All models are fitted on all the available data except for the last year, which is held out for evaluating the $R^2$, so that the evaluation covers a full annual cycle. We call \emph{temperature gain} the increase in held-out $R^2$ from adding (ii) to (i), and \emph{memory gain} the further increase from adding (iii). Since the diagnostic captures only linear and additive effects, the reported gains should be read as lower bounds on the information provided by temperature.

SigAdapt's performance on the forecasting tasks is in line with the respective temperature gains. Temperature adds 0.12 - 0.28 to the $R^2$ in all GFC12 and GFC17 zones except GFC12 zone 9 and GFC17 zones ME and VT (0.09), 0.07 in GFC14, and at most 0.04 in the ENTSO-E countries, whose loads are less temperature-driven since these countries rely mainly on energy sources other than electricity for space heating.\footnote{In 2024, electricity accounted for a minor share of household space heating in all countries of the dataset; see Eurostat, \emph{Energy consumption in households}, \url{https://ec.europa.eu/eurostat/statistics-explained/index.php?title=Energy_consumption_in_households}.}
SigAdapt achieves the best or tied-best RMSE in every series where temperature adds more than 0.1, and all the series where it is outperformed have a temperature gain below 0.1. Conversely, it also outperforms the competitors in two series with a low temperature gain, GFC14 and HU.
These results indicate that SigAdapt outperforms the competitors in settings where the covariate carries substantial information beyond the calendar. The temperature path is moreover informative beyond its current value, as exponentially smoothed temperatures explain up to 21\% of the load variance left unexplained by the instantaneous temperature, most markedly in GFC17.

\begin{table}[ht!]
\centering
\caption{Test performance on the original scale. RMSE and MAE in load
units, MAPE in percentage (\%).
Best model per row in bold.}
\label{tab:results_raw}
\resizebox{\textwidth}{!}{%
\begin{tabular}{ll ccc | ccc | ccc}
\toprule
 & & \multicolumn{3}{c|}{Chronos-2} & \multicolumn{3}{c|}{TiDE} & \multicolumn{3}{c}{SigAdapt} \\
\cmidrule(lr){3-5}\cmidrule(lr){6-8}\cmidrule(lr){9-11}
Dataset & Zone & RMSE & MAE & MAPE & RMSE & MAE & MAPE & RMSE & MAE & MAPE \\
 \midrule
\multirow{11}{*}{GFC12}
 & 1  & 2\,807  & 2\,009  & \textbf{9.85} & 3\,305  & 2\,332  & 11.62 & \textbf{2\,784}  & \textbf{2\,000}  & 10.08 \\
 & 2  & 14\,453 & 10\,347 & 5.73 & 14\,345 & 10\,496 & 5.85 & \textbf{12\,553} & \textbf{8\,838}  & \textbf{4.93} \\
 & 3  & 16\,533 & 12\,379 & 6.32 & 17\,960 & 13\,168 & 6.77 & \textbf{15\,933} & \textbf{11\,839} & \textbf{6.15} \\
 & 4  & 53      & \textbf{38}      & --- & 59      & 42      & --- & \textbf{52}    & \textbf{38}    & --- \\ 
 & 5  & 1\,009  & 736     & 9.48 & 1\,070  & 790     & 10.44 & \textbf{925}   & \textbf{686}   & \textbf{9.17} \\
 & 6  & 16\,313 & 11\,880 & 6.26 & 17\,207 & 12\,663 & 6.73 & \textbf{14\,769} & \textbf{10\,661} & \textbf{5.70} \\
 & 7  & 15\,082 & 11\,072 & 5.67 & 15\,780 & 11\,449 & 5.91 & \textbf{13\,386} & \textbf{9\,779}  & \textbf{5.03} \\
 & 8  & 469     & 349     & 8.72 & 499     & 367     & 9.18 & \textbf{435}   & \textbf{330}   & \textbf{8.32} \\
 & 9  & \textbf{16\,957} & \textbf{11\,449} & --- & 20\,472 & 14\,627 & --- & 17\,827 & 12\,359 & --- \\
 & 10 & 2\,748  & 2\,072  & 7.44 & 3\,084  & 2\,229  & 8.10 & \textbf{2\,606} & \textbf{1\,903} & \textbf{6.93} \\
 & 11 & 12\,876 & 9\,009  & 7.65 & 12\,433 & 9\,161  & 7.91 & \textbf{11\,126} & \textbf{7\,974}  & \textbf{6.78} \\
\midrule
GFC14 & --- & 146 & 98 & 2.91 & 140 & 97 & 2.88 & \textbf{135} & \textbf{93} & \textbf{2.76} \\
\midrule
\multirow{8}{*}{GFC17}
 & CT         & 216 & 153 & 4.41 & 207 & 154 & 4.47 & \textbf{199} & \textbf{142} & \textbf{4.17} \\
 & ME         & \textbf{68}  & \textbf{45}  & \textbf{3.60} & 71 & 49 & 3.93 & 69  & 46  & 3.65 \\
 & NH         & 82  & 58  & 4.39 & 80 & 59 & 4.49 & \textbf{67}  & \textbf{50}  & \textbf{3.78} \\
 & RI         & 66  & 46  & 4.83 & 62 & 44 & 4.67 & \textbf{54}  & \textbf{38}  & \textbf{4.13} \\
 & VT         & \textbf{39}  & \textbf{27}  & \textbf{4.54} & 42 & 31 & 5.18 & 42  & 31  & 5.16 \\
 & SEMASS     & 132 & 95  & 5.67 & 125 & 88 & 5.32 & \textbf{118} & \textbf{85}  & \textbf{5.18} \\
 & WCMASS     & 122 & 89  & 4.68 & \textbf{114} & \textbf{84} & \textbf{4.45} & \textbf{114} & \textbf{84}  & 4.51 \\
 & NEMASSBOST & 198 & 140 & \textbf{4.84} & 209 & 152 & 5.29 & \textbf{186} & \textbf{138} & 4.87 \\
 \midrule
\multirow{6}{*}{ENTSO-E}
 & AT & \textbf{1\,514}  & \textbf{897}   & \textbf{3.21}  & 1\,558  & 952    & 3.41  & 1\,525  & 1\,009 & 3.56 \\
 & BE & \textbf{1\,673}  & \textbf{1\,157} & \textbf{3.04}  & 1\,740  & 1\,159 & \textbf{3.04}  & 1\,786  & 1\,247 & 3.29 \\ 
 & DE & 13\,603 & 7\,649 & 3.58  & \textbf{10\,847} & \textbf{6\,703} & \textbf{3.13}  & 12\,208 & 7\,781 & 3.64 \\
 & HU & 1\,086  & 645   & 3.35  & 1\,048  & \textbf{626}   & \textbf{3.28}  & \textbf{987}    & 676    & 3.49 \\
 & LU & \textbf{294}    & \textbf{209}   & \textbf{13.02} & 314    & 232    & 14.31 & 331    & 252    & 15.99 \\
 & NL & 2\,049  & 1\,296 & 2.51  & \textbf{2\,001}  & \textbf{1\,285} & \textbf{2.50}  & 2\,200  & 1\,458 & 2.80 \\
 \bottomrule
\end{tabular}%
}
\end{table}

\subsection{Results on the electricity price dataset}
\label{app:results_epf}

Results on the original scale are reported in Table~\ref{tab:epf_orig}. SigAdapt outperforms TiDE but is outperformed by Chronos-2 in all markets except PJM, and DE in terms of RMSE.
Its errors are concentrated around price spikes, which a linear model on signature features captures less well than load dynamics. Improving the treatment of such jumps is a direction for future work.

\begin{table}[ht!]
\centering
\caption{EPF test performance on the original scale.  Best model per row in bold.}
\label{tab:epf_orig}
\resizebox{0.52\textwidth}{!}{%
\begin{tabular}{l cc | cc | cc}
\toprule
Models & \multicolumn{2}{c|}{Chronos-2} & \multicolumn{2}{c|}{TiDE} & \multicolumn{2}{c}{SigAdapt} \\
\cmidrule(lr){2-3}\cmidrule(lr){4-5}\cmidrule(lr){6-7}
Market & RMSE & MAE & RMSE & MAE & RMSE & MAE \\
\midrule
FR  & \textbf{10.41} & \textbf{4.06} & 13.01 & 5.71 & 11.03 & 4.34 \\
DE  & 6.21 & \textbf{3.43} & 7.86 & 4.59 & \textbf{6.16} & 3.57 \\
BE  & \textbf{15.22} & \textbf{6.05} & 17.55 & 7.74 & 15.50 & 6.36 \\
NP  & \textbf{3.65} & \textbf{1.86} & 3.90 & 2.08 & 3.66 & 2.07 \\
PJM & 5.39 & 3.18 & 5.87 & 3.49 & \textbf{5.33} & \textbf{3.05} \\
\bottomrule
\end{tabular}%
}
\end{table}

\end{document}